\documentclass[11pt, oneside, a4paper, reqno]{amsart}
\usepackage[scale = .75]{geometry}

\usepackage[utf8]{inputenc}

\usepackage{amsmath, amssymb}
\usepackage{stmaryrd}
\usepackage{mathtools}
\usepackage{tikz-cd}
\usepackage{enumitem}
\usepackage[ruled]{algorithm2e}

\usepackage[hidelinks]{hyperref}

\swapnumbers
\newtheorem{prop}{Proposition}[section]
\newtheorem{thm}[prop]{Theorem}
\newtheorem{cor}[prop]{Corollary}
\newtheorem{lemma}[prop]{Lemma}
\newtheorem*{thm*}{Theorem}

\theoremstyle{definition}

\newtheorem{example}[prop]{Example}

\numberwithin{equation}{section}

\swapnumbers % AMS number swapping overrides all theorem styles
\newtheoremstyle{stepstyle}%
{}{}% space above and below
{}{}% body font and indent
{}{.}% theorem head font and punctuation
{ }% space after theorem head
{\thmnumber{#2}.\ \bfseries\boldmath\thmnote{#3}}

\theoremstyle{stepstyle}
\newtheorem{step}[prop]{}

\newtheoremstyle{examplestyle}%
{}{}% space above and below
{}{}% body font and indent
{\bfseries\boldmath}{.}% theorem head font and punctuation
{ }% space after theorem head
{\thmname{#1}\thmnumber{ #2}:\ \thmnote{#3}}

\theoremstyle{examplestyle}

\newcounter{halfsection}[section]
\renewcommand{\thehalfsection}{Example \Alph{halfsection}}
\newcommand\halfsection[1]{\refstepcounter{halfsection}%
\section*{\thehalfsection. #1}%
}

\def\diagramscale{.45}
\tikzset{%
vcenter/.style= {baseline = {([yshift = -.5ex, scale=\diagramscale](0, #1/2+.5))}},%
functor fill/.style = {fill = black!30!white, opacity = .6},%
functor draw/.style = {line width = 2pt, draw = black!50!white, opacity = .6},%
}
\newcommand{\tikzgrid}[1]{%
\foreach \x in {0, ..., #1} {
\draw[dotted] (\x+.5, .5) -- (\x+.5, #1+.5);
\draw[dotted] (.5, \x+.5) -- (#1+.5, \x+.5);
}
}

\def\emph{\textbf}
\def\parameterization{\mathcal{T}}

\DeclareMathOperator{\coker}{coker}
\DeclareMathOperator{\rank}{rank}
\def\Fun{\textnormal{Fun}}
\def\id{\textnormal{id}}
\def\Nat{\textnormal{Nat}}
\def\op{{\textnormal{op}}}
\def\supp{\textnormal{supp}}
\def\th{{\textnormal{th}}}
\def\vect{\textnormal{vect}_K}

\def\Up{\textnormal{Up}}

\DeclareMathOperator{\Max}{Max}
\DeclareMathOperator{\Min}{Min}

\def\leqJ{\preccurlyeq}
\def\ltJ{\prec}
\def\geqJ{\succcurlyeq}

\title[Koszul complexes and relative homological algebra]{Koszul complexes and relative homological algebra of functors over posets}

\author[W. Chach\'olski]{Wojciech Chach\'olski}
\address{Mathematics, KTH, S-10044 Stockholm, Sweden}
\email{Wojciech Chach\'olski <wojtek@kth.se>}
\email{Andrea Guidolin <guidolin@kth.se>}
\email{Isaac Ren <isaacren@kth.se>}
\email{Martina Scolamiero <scola@kth.se>}
\address{Max Planck Institute for Mathematics in the Sciences, 04103 Leipzig, Germany}
\email{Francesca Tombari <francesca.tombari@mis.mpg.de>}
\thanks{Communicated by Ulrich Bauer}
\thanks{Partially supported by VR, the Wallenberg AI, Autonomous System and Software Program (WASP)
and Data Driven Life Science (DDLS) program funded by Knut and Alice Wallenberg Foundation,  MultipleMS funded by the European Union under the Horizon 2020 program, grant agreement 733161,  dBRAIN collaborative project at digital futures at KTH}

\author[A. Guidolin]{Andrea Guidolin}
\author[I. Ren]{Isaac Ren}
\author[M. Scolamiero]{Martina Scolamiero}
\author[F. Tombari]{Francesca Tombari}


\begin{document}

\begin{abstract}
Under certain conditions, Koszul complexes can be used to calculate relative Betti diagrams of vector space-valued functors indexed by
a poset, without the explicit computation of global minimal relative resolutions. 
In relative homological algebra of such functors, free functors are replaced by an arbitrary family of functors.
Relative Betti diagrams encode the multiplicities of these functors in minimal relative resolutions.
In this article we provide conditions under which grading the chosen family of functors leads to explicit Koszul complexes whose homology dimensions are the relative Betti diagrams, thus giving a scheme for the computation of these numerical descriptors.
\end{abstract}

\maketitle

\section{Introduction}

\subsection{Relative homological algebra and Koszul complexes}
Recently, there has been intense research activity~\cite{AENY2019, BBH2022,BOO2021} relating topological data analysis (TDA) with relative homological algebra \cite{ARS1995, auslander1993relative, eilenberg1965foundations, enochs2000relative, Hochschild1956} of finite dimensional vector space representations of finite posets. 
One way of thinking about such representations is as modules over certain finite dimensional algebras, allowing to phrase homological properties of modules in terms of homological properties of the algebra. 
This is also a common way for expressing relative homological properties of representations.
In this case, endomorphism algebras of some chosen modules are often utilized (see e.g.\ \cite{ARS1995,BBH2022}). 
We however take another perspective, central in TDA, and study representations as functors. 
Local Koszul complexes \cite{tameness} (see also \cite[Sect.~6]{matsumura1989commutative} and \cite[Sect.~1]{miller2005combinatorial}) are the key reasons why viewing representations as functors is particularly convenient in TDA. 
In favourable situations, these complexes can be used for calculations of Betti diagrams without explicitly constructing global minimal resolutions, especially if we have some control over the sets of parents of elements in the indexing poset.

The aim of this article is to show that Koszul complexes can also be used to calculate relative Betti diagrams.
To do this, we translate from a relative homological algebra of functors indexed by a given poset to the standard homological algebra of functors indexed by a new poset.
We are interested in translations for which the homology of the Koszul complexes over the  new poset computes the relative Betti diagrams in a way that can be implemented in software (see for example \cite{lowerhooks}).
One of our main contributions in this article is a translation construction which allows us to characterize the set of poset elements (called the degeneracy locus) where local Koszul complexes fail to calculate the relative Betti diagrams. 
Importantly, in all examples that are currently relevant for TDA the degeneracy loci are negligible. 

Aiming at a self-contained presentation, and believing that the tools we introduce may be of interest to a broader audience, we favor an explicit classical formulation of relative homological algebra~\cite{eilenberg1965foundations}.

\subsection{Relative homological approximations}
One way of understanding an object is by approximating it.
In this article we focus on approximations coming from relative homological algebra (see~\cite{Hochschild1956}), where arbitrary objects in an abelian category are approximated by finite direct sums of objects of a chosen collection $\mathcal{P}$.
For example, one may ask whether a given object $M$ is isomorphic to a finite direct sum of elements in $\mathcal{P}$.
To answer this question we look at the category $\bigoplus\mathcal{P}\downarrow_M$ of morphisms from finite direct sums of elements in $\mathcal{P}$ into $M$. 
We are interested in situations where this category has a certain distinguished object $C_0 \to M$: for example, a terminal object.
However, cases where the category $\bigoplus \mathcal{P} \downarrow_M$ has a terminal object are not that common. More common is the existence of $C_0 \to M$ in the category $\bigoplus\mathcal{P}\downarrow_M$ satisfying the following two conditions.
The first condition requires that every object in $\bigoplus \mathcal{P} \downarrow_M$ maps to $C_0\to M$, in which case $C_0 \to M$ is called a \emph{$\mathcal{P}$-epimorphism} (see~\ref{Rel:Exactness}).
This of course holds if $C_0 \to M$ is terminal in $\bigoplus \mathcal{P} \downarrow_M$, where not only the existence but also  the uniqueness of morphisms into  $C_0 \to M$ is required.
The second condition substitutes this uniqueness by a minimality condition: every endomorphism of $C_0 \to M$ is an isomorphism.
A morphism $C_0 \to M$ satisfying these two conditions (being a minimal $\mathcal{P}$-epimorphism) is called a \emph{minimal $\mathcal{P}$-cover} of $M$ (see~\ref{Rel:Covers}).
If it exists, then the minimal $\mathcal{P}$-cover of $M$ is unique up to isomorphism. 
We regard the minimal $\mathcal{P}$-cover of $M$ as an approximation of $M$ by finite direct sums of elements in $\mathcal{P}$.
In particular, $M$ is isomorphic to such a direct sum if, and only if, its  minimal $\mathcal{P}$-cover is an isomorphism.

If finite direct sums of elements in $\mathcal{P}$ are uniquely determined by the multiplicities of their summands, then we say that $\mathcal{P}$ is \emph{independent} (see~\ref{Rel:Freeness}).
In this case, there is a unique function $\beta_{\mathcal{P}}^{0}M\colon \mathcal{P}\to\mathbb{N}$ whose value at $A$ in $\mathcal{P}$ is the multiplicity of $A$ in the minimal $\mathcal{P}$-cover $C_0$ of $M$. 
The function $\beta_{\mathcal{P}}^{0}M$ is called the \emph{$0^{\th}$ $\mathcal{P}$-Betti diagram} of $M$ and is a numerical descriptor of the approximation of $M$ by finite direct sums of elements in $\mathcal{P}$.
This is only the $0^\th$ step.
We can continue by considering the minimal $\mathcal{P}$-cover of the kernel of $C_0\to M$.
By doing this inductively we obtain a chain complex $\cdots \to C_n \to \cdots \to C_0 \to M$ called a \emph{minimal $\mathcal{P}$-resolution} of $M$. 
Each object $C_d$ in this complex is described by a function $\beta^d_{\mathcal{P}} M \colon \mathcal{P} \to \mathbb{N}$ called the \emph{$d^\th$ $\mathcal{P}$-Betti diagram of $M$}, where, for all $A$ in $\mathcal{P}$, the number $\beta^d_{\mathcal{P}} M(A)$ is the multiplicity of $A$ in $C_d$ (see~\ref{Rel:Betti}).

Thus every independent collection $\mathcal{P}$, for which all objects admit minimal $\mathcal{P}$-covers, leads to numerical invariants in the form of $\mathcal{P}$-Betti diagrams.
The focus of this article is to provide a method of constructing such collections $\mathcal{P}$ and to discuss a strategy to effectively compute the associated Betti diagrams in the category $\Fun(I,\vect)$ of functors indexed by a finite poset $I$ with values in finite dimensional $K$-vector spaces.

\subsection{Related work}
\label{sec:relatedwork}
Since vector space valued functors indexed by posets have played an important role across many areas of mathematics, their categories are well studied, and in particular their homological properties are well understood: see for example~\cite{ARS1995} and references therein, where such functors are interpreted as modules over finite dimensional algebras.
Recently there has also been a growing interest in independent collections in $\Fun(I,\vect)$ \cite{BOO2021,AENY2023,BBH2022} and associated numerical invariants coming from the TDA community, where such functors are commonly called persistence modules, and the poset $I$ is usually $\mathbb{R}^n$ with the product order. 
In particular, there are several 
publications and implementations of various algorithms to calculate numerical homological invariants of such functors (e.g.\ classical Betti diagrams, rank functions, generalized Euler characteristics) and using them in applications, see for example~\cite{LW2015, LW2022, LSCBO2023, MP2023, OS2024, RD2020}.

The idea of approximating persistence modules by objects in a chosen collection appears in \cite{kim2021generalized,AENY2019,BOO2021,AENY2023,BBH2022}.
In \cite{kim2021generalized}, without relying on homological algebra, approximations in terms of the intervals in the poset are shown to yield a generalization of rank functions, an invariant of persistence modules. In \cite{AENY2019}, the authors approximate persistence modules in a Grothendieck group generated by interval modules (called spread modules in \cite{BBH2022} and in our work, see \ref{S:SpreadModules}) using quiver representation theory without resorting to homological algebra methods. 
Decompositions of rank functions into linear combinations of rank functions of simpler modules, such as rectangles and lower hooks, are explained in~\cite{BOO2021}. Importantly, one such decomposition (called minimal rank decomposition) is defined via relative homological algebra of persistence modules with respect to the collection of lower hooks, which we consider in \ref{E:LowerHooks}.
In \cite{AENY2023}, relative homological algebra with respect to interval modules is studied using the language of representations of finite dimensional algebras. 
For homological algebra results for persistence modules over more general posets, we can also look at \cite{BM2021} and \cite{Miller2020}. In the former, the authors prove various standard homological algebra results for persistence modules, with a focus on the one-parameter case (where the poset is the real line).
In the latter, the author studies presentations of persistence modules, with focus on infinite posets, and concludes with a Hilbert syzygy-type theorem. 
The authors of \cite{BBH2022} present
a way of obtaining relative projective resolutions from standard projective resolutions by considering modules over the ring $\text{End}(\bigoplus \mathcal{P})^\op$ and finding a resolution of $\hom(\bigoplus \mathcal{P}, M)$.
This gives an automatic correspondence of exact structures, but there is no reason for the standard projective resolution of $\hom(\bigoplus \mathcal{P}, M)$ to be more easily computed than the relative projective resolution of $M$.

Our original motivation was to find a homological interpretation of the approximations presented in these articles that would allow for explicit calculations. This article builds on the realization that local Koszul complexes can be used for this purpose.
Koszul complexes of modules are well studied, see e.g.~\cite[Sect.~6]{matsumura1989commutative}. 
Koszul complexes of graded modules, as presented in \cite{miller2005combinatorial}, have a direct interpretation in terms of persistence modules, based on the equivalence between the category of $n$-graded modules over the polynomial ring $S\coloneqq K [x_1,\ldots ,x_n]$ and the category $\Fun(\mathbb{N}^n,\vect)$, where $\mathbb{N}^n$ is endowed with the product order.
Here we sketch this connection, referring to \cite[Sect.~3]{guidolin2023morse} for more details. 
The Koszul complex $\mathcal{K}$ as defined for example in \cite[Def.~1.26]{miller2005combinatorial} is a minimal free resolution of the field $K$ in the category of $n$-graded $S$-modules \cite[Prop.~1.28]{miller2005combinatorial}.
Given an $n$-graded $S$-module $M$, the Koszul complex of $M$ at  $a\in\mathbb{N}^n$ is the grade $a$ part of the $n$-graded chain complex $M \otimes_S \mathcal{K}$. 
Its $d^\th$ homology module is the grade $a$ part of the $n$-graded $S$-module $\operatorname{Tor}_i^{S}(M,K)$, whose dimension is the $d^\th$ Betti diagram (with respect to the standard projectives) of $M$ at $a$.
This value of the Betti diagram can be computed, for a grade $a\in \mathbb{N}^n$, by looking only at $M$ restricted to the grades in a small subposet of $\mathbb{N}^n$ isomorphic to $\{0<1\}^k$ for some $k\le n$, whose elements are $a$ together with its parents and all their meets.
In this work, we study local Koszul complexes of functors indexed by more general finite posets (see \ref{Std:KoszulComplexes}), more details on which can be found in \cite{tameness}, and use them in the context of relative homological algebra. 

The above approach starts with a global resolution of $K$ and considers local parts of $M \otimes_S \mathcal{K}$ at different grades. Such a global resolution however may not have an explicit formula, which typically would require global regularity assumptions on the indexing poset. Requiring local regularity only at some  grades is a much weaker assumption. 
In such grades we have an explicit Koszul complex construction
which we exploit. For that purpose,
the functorial perspective is more advantageous than  the language of representations of finite dimensional algebras used in \cite{BOO2021,AENY2023,BBH2022}. Our aim is to use Koszul complexes not only as a theoretical tool, but as a computable construction to determine relative Betti diagrams.

\subsection{Koszul complexes and thinness}
Homological algebra characterizes relevant objects by  universal properties, regardless of the category to which they belong, instead of by concrete constructions which may vary greatly with the category.
This offers a powerful language in which to express  and solve mathematical problems.
However, extracting calculable invariants, which is essential in TDA, from universally defined objects is often difficult, costly, and sometimes not even feasible.
Therefore the challenge is to introduce a version of homological algebra which retains some conceptual power while allowing for some of the relevant invariants to be described explicitly and algorithmically.

We strike a balance between mathematical and computational effectiveness by using a strategy built on the following two observations.
The first is related to the standard independent collection $\mathcal{S} = \{K(a,-) \mid a \in I\}$ consisting of all free functors on one generator (see beginning of Section~\ref{S:StdHomAlg}).
This collection is always independent and, when $I$ is finite, every functor in $\Fun(I,\vect)$ admits a minimal $\mathcal{S}$-resolution.
The first key observation is that under additional assumptions, for example that $I$ is an upper semilattice, for every element $a$ in $I$ and every functor $F \colon I \to \vect$, there exists an explicit chain complex $\mathcal{K}_a F$, called \textbf{Koszul complex} (see~\ref{Std:KoszulComplexes} and discussion in \ref{sec:relatedwork}), whose homology in degree $d$ has dimension equal to the value of the Betti diagram $\beta^d_{\mathcal{S}} F(K(a,-))$ (see Theorem~\ref{T:BettiKoszul} and \cite[Section 10]{tameness}).
Because of its explicitness, the Koszul complex is an effective tool for computing the values of the Betti diagram with respect to the collection $\mathcal{S}$.
 
To have an analogous construction for a more general collection $\mathcal{P}$ in $\Fun(I,\vect)$, we require a grading of its elements.
Thus, instead of an unstructured collection, our starting point is a functor $\parameterization \colon J^\op \to \Fun(I, \vect)$ indexed by a finite poset $(J, \leqJ)$.
Such a functor leads to a collection $\mathcal{P} := \{\parameterization(a) \mid a \in J,\ \parameterization(a) \neq 0\}$.
This grading allows for a translation between functors indexed by $I$ and functors indexed by $J$ via a pair of standard adjoint functors (see~\ref{D:Adjunction} and~\ref{sadgsghj}):
\[
\begin{tikzcd}
\Fun(I,\vect)
 \ar[bend left=10]{r}{\mathcal R}
&
\Fun(J,\vect)
 \ar[bend left=10]{l}{\mathcal L}
\end{tikzcd}
\]
where $\mathcal R := \Nat_I(\parameterization, -)$ assigns to $M \colon I \to \vect$ the functor $\Nat_I(\parameterization(-), M) : J \to \vect$. The functor $\parameterization$ is called \textbf{thin} if the unit of this adjunction $\eta_a \colon K(a,-) \to \mathcal{R} \mathcal{L} K(a,-)$ is an epimorphism for every $a$ in $J$.  
The second key observation is recognizing the importance of the thinness assumption. 
This assumption enables us to translate between the homological properties  
of $M$ in  $\Fun(I,\vect)$ relative to $\mathcal{P}$ and the standard homological properties of $\mathcal{R} M$ in $\Fun(J,\vect)$.

\subsection{Results} 
Let $I$ and $J$ be finite posets, and
$\parameterization \colon J^\op \to \Fun(I, \vect)$ a thin functor. 
We give a self-contained proof that the collection $\mathcal{P}=\{\parameterization(a)\mid a\in J,\  \parameterization(a)\neq 0\}$ is independent (\ref{P:ThinIndependence}). 
Moreover, we prove that the adjunction between $\mathcal{L}$ and $\mathcal{R}$ leads to a one-to-one correspondence between minimal covers: 
for every functor $M$ in $\Fun(I, \vect)$, a natural transformation $C_0 \to \mathcal{R} M$ is a minimal cover in $\Fun(J,\vect)$ if, and only if, its adjoint $\mathcal{L} C_0 \to M$ is a minimal $\mathcal{P}$-cover in $\Fun(I, \vect)$ \textnormal{(Theorem \ref{T:MinCovers})}. 
Although we show that every functor $M$ in $\Fun(I, \vect)$ admits a minimal $\mathcal{P}$-resolution \textnormal{(Corollary \ref{P:ThinAcyclicity}(1))}, this resolution cannot in general be obtained via the adjunction. This means that in general the $\mathcal{P}$-Betti diagrams of $M$ differ from the standard Betti diagrams of $\mathcal{R} M$.
Elements of $J$
where this happens are called degenerate and form the degeneracy locus (\ref{D:DegeneracyLocus}). If $J$ is an upper semilattice, 
our main result (Corollary \ref{C:BettiFromDegeneracy}) identifies elements outside 
the degeneracy locus, for which we can use Koszul complexes to calculate the values of $\mathcal{P}$-Betti diagrams. From a computational perspective, we obtain Algorithm \ref{A:lgo} for computing relative Betti diagrams, given a thin functor.

\begin{algorithm}[ht]
\label{A:lgo}
\caption{Betti diagrams relative to a parameterization $\mathcal{T}$}
\KwIn{$I$ finite poset, $J$ finite upper semilattice, \\
$\parameterization \colon J^{\op} \to \Fun(I, \vect)$ satisfying the conditions of Corollary \ref{C:BettiFromDegeneracy}, \\
$M \colon I \to \vect$ functor, \\
$a \in J$ poset element such that $\parameterization(a) \neq 0$.}
\KwOut{the relative Betti diagrams of $M$ at $\parameterization(a)$ relative to $\{\parameterization(b) \mid b \in J, \; \parameterization(b) \neq 0\}$.}

$\mathcal{U}(a) \leftarrow \{a_1, \ldots, a_n\}$ parents of $a$\tcp*{parents defined in \ref{Std:KoszulComplexes}}
$J_0 \leftarrow \{\varnothing\}$\;
\lFor{$d = 1, \ldots, n$}
{$J_d \leftarrow \{S \subseteq \mathcal{U}(a) \mid S \textnormal{ bounded below}, \; |S| = d\}$}
\lFor(\tcp*[f]{$\bigwedge \varnothing = a$ by convention}){$S \in \bigcup_{d = 0}^n J_d$}
{$\mathcal{R}M(\bigwedge S) \leftarrow \Nat(\mathcal{T}(\bigwedge S), M)$}
\For{$d = 0, \ldots, n - 1$}{%
\For{$S \in J_d$, $S' \in J_{d + 1}$, and $S \subset S'$}{%
Write $i_0 < \cdots < i_d$ such that $S' = \{a_{i_0}, \ldots, a_{i_d}\}$\;
$k_{S, S'} \leftarrow k$ such that $S' \setminus S = a_{i_k}$\;
$b \leftarrow \bigwedge S$; $c \leftarrow \bigwedge S'$\;
$\partial_{S, S'} \leftarrow \big(\Nat(\mathcal{T}(b > c), M) \colon \mathcal{R}M (c) \to \mathcal{R}M (b)\big)$\;
}
}
\tcp{Construction of chain complex $(\mathcal{K}, \partial)$ from \ref{Std:KoszulComplexes}}
$\mathcal{K} \leftarrow \left( \bigoplus_{S \in J_d} \mathcal{R}M(\bigwedge S) \right)_{d = 0, \ldots, n}$\;
$\partial \leftarrow \left(\sum_{S \in J_d, \; S' \in J_{d + 1}} (-1)^{k_{S, S'}} \partial_{S, S'} \right)_{d = 0, \ldots, n - 1}$\;
$\beta \leftarrow \dim H_*(\mathcal{K}, \partial)$\;
\Return{$\beta$}
\end{algorithm}

Furthermore, we  produce examples of thin functors $\parameterization \colon J^\op \to \Fun(I, \vect)$, and we apply our results to explicitly calculate the associated relative Betti diagrams.
We observe how numerical invariants introduced in other works, such as the rank invariant over lower hooks of \cite{BOO2021} (\ref{E:LowerHooks}) and the single-source homological spread invariant of \cite{BBH2022} (\ref{E:SingleSourceSpreadModules}), are included in our general framework.

\subsection{Organization}
We begin in Section~\ref{S:RelHomAlg} by introducing the theory of relative homological algebra, including the notions of independence and acyclicity, covers and resolutions. 
Next, in Section~\ref{S:StdHomAlg} we recall standard homological algebra for functors indexed by finite posets. In particular, when the poset is an upper semilattice, we show how to compute Betti diagrams explicitly via Koszul complexes.
In Section~\ref{S:FiltrationsSubfunctors}, we discuss in more depth resolutions and Betti diagrams of certain functors called filtrations. In particular, we study subfunctors, and discuss several important notions: supports of functors, upsets, and antichains.

Section~\ref{S:GradedRelHomAlg} is the central part of this work, where we introduce gradings of the form $\parameterization \colon J^\op \to \Fun(I, \vect)$. We define the key properties of thinness and flatness, which will allow us to compute relative Betti diagrams over $I$ from standard Betti diagrams over $J$.

Finally, we connect to other works in Section~\ref{S:Examples} by exhibiting how rectangles, lower hooks, and spread modules, among others, are handled in our framework.

\subsection*{Acknowledgements}
We would like to thank Luis Scoccola for valuable discussions and the anonymous referees for their extensive reports and constructive suggestions.

\section{Relative homological algebra}
\label{S:RelHomAlg}

In this section we recall a setup for relative homological algebra in an abelian category $\mathcal M$.
This subject, initiated by Hochschild \cite{Hochschild1956}, has been extensively developed in the context of the representation theory of finite dimensional algebras (see for example~\cite{ARS1995, auslander1993relative}). The explicit and self-contained treatment we provide here has the advantage of not assuming familiarity of the reader with such theory.
We fix a collection $\mathcal{P}$ of objects in $\mathcal M$.
What follows are fundamental notions of homological algebra relative to $\mathcal{P}$.

\begin{step}[Freeness]
\label{Rel:Freeness}
An object in $\mathcal M$ is called \emph{$\mathcal{P}$-free} if it is isomorphic to a finite direct sum of elements in $\mathcal{P}$. For example, the zero object is $\mathcal{P}$-free.
The collection $\mathcal{P}$ is called \emph{independent} if, for every $\mathcal{P}$-free object $F$, there is a unique function $\beta \colon \mathcal{P} \to \mathbb{N}$ whose \emph{support} $\supp \, \beta:=\{A \mid \beta(A) \neq 0\}$ is finite and 
for which $F$ is isomorphic to $\bigoplus_{A\in \mathcal{P}} A^{\beta(A)}$.
If $\mathcal{P}$ is independent, then all of its elements are nonzero, and if two of its elements are isomorphic, then they must be equal.
For example, if $\mathcal M$ is a Krull-Schmidt category and $\mathcal{P}$ is a collection of indecomposables, then it is independent if and only if
no two elements in $\mathcal{P}$ are isomorphic.
\end{step}

\begin{step}[Exactness]\label{Rel:Exactness}
Two composable morphisms $f \colon M \to N$ and $g \colon N \to L$ in $\mathcal{M}$ are said to form a \emph{$\mathcal P$-exact sequence} if the following sequence of abelian groups is exact for every $A$ in $\mathcal{P}$:
\[
\begin{tikzcd}[column sep=40pt]
\hom(A, M) \ar{r}{\hom(A, f)}
& \hom(A, N) \ar{r}{\hom(A, g)}
& \hom(A, L).
\end{tikzcd}
\]
A sequence of composable morphisms $\cdots \to M_d \to M_{d-1} \to \cdots$ is called \emph{$\mathcal P$-exact} if every two of its consecutive morphisms form a $\mathcal P$-exact sequence.

If the sequence $M \xrightarrow{f} N \to 0$ is $\mathcal P$-exact, then $f$ is called a \emph{$\mathcal{P}$-epimorphism}. A morphism $f \colon M \to N$ is a $\mathcal{P}$-epimorphism if, and only if, $\hom(A, f) \colon \hom(A, M) \to \hom(A, N)$ is surjective for every $A$ in $\mathcal{P}$ (i.e.\ every $g \colon A \to N$ can be expressed as the composition of some morphism $A \to M$ with $f$).
\end{step}

\begin{step}[Projectives]
\label{Rel:Projectives}
An object $C_0$ in $\mathcal M$ is called \emph{$\mathcal{P}$-projective} if, for every 
$f \colon M \to N$ and $g \colon N \to L$ that form a 
$\mathcal{P}$-exact sequence, the following is an exact sequence of abelian groups:
\[
\begin{tikzcd}[column sep=40pt]
\hom(C_0, M) \ar{r}{\hom(
C_0, f)}
& \hom(C_0, N) \ar{r}{\hom(C_0, g)}
& \hom(C_0, L).
\end{tikzcd}
\]
Note that if $M$ is $\mathcal{P}$-projective, then $f \colon M \to N$ and $g \colon N \to L$ form a $\mathcal P$-exact sequence if, and only if, the composition $gf$ is null and the induced morphism $M \to \ker(g)$ is a $\mathcal P$-epimorphism.
Every $\mathcal{P}$-free object is $\mathcal{P}$-projective. If every $\mathcal{P}$-projective is $\mathcal{P}$-free, then the collection $\mathcal{P}$ is called \emph{acyclic}.
\end{step}

\begin{step}[Covers]
\label{Rel:Covers}
A \emph{$\mathcal{P}$-cover} of an object $M$ is a $\mathcal{P}$-epimorphism $C_0 \to M$ where $C_0$ is $\mathcal{P}$-projective.
We say that the category $\mathcal{M}$ \emph{has enough $\mathcal{P}$-projectives} if every object has a $\mathcal{P}$-cover.

A morphism between a $\mathcal{P}$-cover $C_0 \to M$ and a $\mathcal{P}$-cover $D_0\to M$ is a morphism $f\colon C_0\to D_0$ in $\mathcal M$ for which the following diagram commutes:
\[
\begin{tikzcd}[row sep=5pt]
C_0 \ar{dd}[swap]{f} \ar{dr} \\
& M \\
D_0 \ar{ur}
\end{tikzcd}
\]
Such a morphism is an isomorphism if $f$ is an isomorphism in $\mathcal M$.
There is always a morphism between two $\mathcal{P}$-covers $C_0 \to M$ and $D_0 \to M$, since $C_0$ is $\mathcal{P}$-projective and $D_0 \to M$ is a $\mathcal{P}$-epimorphism.

A $\mathcal{P}$-cover $C_0\to M$ is called \emph{minimal} if all its endomorphisms are isomorphisms. Every two minimal $\mathcal{P}$-covers of $M$ are isomorphic: if $C_0 \to M \leftarrow D_0$ are two minimal $\mathcal{P}$-covers, then in any sequence  $C_0 \to D_0 \to C_0 \to D_0$ of morphisms of the $\mathcal{P}$-covers, the first and second ones compose into an isomorphism by minimality, as well as the second and third ones, which implies that all these morphisms are isomorphisms. 
\end{step}

\begin{step}[Resolutions]
\label{Rel:Resolutions}
A sequence of composable morphisms in $\mathcal M$ of the form $\cdots \to C_1 \to C_0 \to M\to 0$ is called a \emph{$\mathcal P$-resolution} of $M$ if it is $\mathcal P$-exact and, for every $d\geq 0$, the object $C_d$ is $\mathcal P$-projective.
Equivalently, the sequence $\cdots \to C_1 \to C_0 \to M \to 0$ is a $\mathcal P$-resolution of $M$ if the following three conditions are satisfied: (a) for every $d\geq 0$, the object $C_d$ is $\mathcal P$-projective; (b) the composition of every two consecutive morphisms is $0$; and (c) for every $d\geq 0$, the induced morphism $C_d\to \ker(C_{d - 1} \to C_{d - 2})$ is a $\mathcal P$-epimorphism, where $C_{-2}=0$ and $C_{-1} = M$.

A $\mathcal P$-resolution of $M$ is also written as a morphism of chain complexes, by $C \to M$, where $M$ denotes a chain complex concentrated in degree $0$ and $C = (\cdots \to C_1 \to C_0)$.
An object $M$ for which there is a $\mathcal P$-resolution $C \to M$ is called \emph{$\mathcal P$-resolvable}. If there are enough $\mathcal P$-projectives, then all objects are $\mathcal P$-resolvable.
 
A morphism from a $\mathcal{P}$-resolution $C \to M$ to a $\mathcal{P}$-resolution $D \to M$ is a sequence of morphisms $f_d \colon C_d \to D_d$ for $d \geq 0$ for which the following diagram commutes:
\[
\begin{tikzcd}[row sep=5pt]
\cdots \ar{r} & C_1 \ar{dd}[swap]{f_1} \ar{r} & C_0 \ar{dd}[swap]{f_0} \ar{dr} \\
& & & M \\
\cdots \ar{r} & D_1 \ar{r} & D_0 \ar{ur}
\end{tikzcd}
\]
Such a morphism is an isomorphism if $f_d$ is an isomorphism for every $d \geq 0$.
As with $\mathcal{P}$-covers, there is always a morphism between every two $\mathcal P$-resolutions of $M$.
 
A $\mathcal P$-resolution $C \to M$ is called \emph{minimal} if all of its endomorphisms are isomorphisms. A $\mathcal P$-resolution $C \to M$ is minimal if, and only if, for every $d \geq 0$, the morphism $C_d \to \ker(C_{d-1} \to C_{d-2})$ is a minimal $\mathcal{P}$-cover, where $C_{-1} = M$ and $C_{-2} = 0$. Thus, if every object has a minimal $\mathcal P$-cover, then every object has a minimal $\mathcal P$-resolution that can be constructed inductively by taking minimal $\mathcal P$-covers of successive kernels.
Every two minimal $\mathcal P$-resolutions of $M$ are isomorphic.
\end{step}

\begin{step}[Betti diagrams]
\label{Rel:Betti}
Suppose $\mathcal{P}$ is independent (see~\ref{Rel:Freeness}) and acyclic (see~\ref{Rel:Projectives}). Let $C \to M$ be a minimal $\mathcal P$-resolution. By acyclicity of $\mathcal{P}$, for every $d \geq 0$, the object $C_d$ is $\mathcal{P}$-free. Thus, by independence of $\mathcal{P}$, there is a unique function, denoted by 
$\beta^d_{\mathcal P} M \colon \mathcal{P} \to \mathbb{N}$ and 
called the \emph{$d^{\th}$ $\mathcal P$-Betti diagram of $M$}, for which 
$C_d$ is isomorphic to
\[
\bigoplus_{A \in \mathcal{P} } A^{\beta^d_{\mathcal P} M(A)}
\]
Since any two minimal $\mathcal P$-resolutions of $M$ are isomorphic, the $\mathcal P$-Betti diagrams of $M$ do not depend on the choice of a minimal $\mathcal{P}$-resolution, and are invariants of the isomorphism type of $M$. One should note however that $\mathcal P$-Betti diagrams are only defined for objects that have minimal $\mathcal P$-resolutions. 

To calculate the values of the $d^{\th}$ $\mathcal P$-Betti diagram of $M$, a standard strategy is to develop tools to calculate the $0^\th$ $\mathcal P$-Betti diagram and then use the fact that, for a minimal resolution $C \to M$, we have the following sequence of equalities:
\[
\beta^d_{\mathcal P}M=\beta^{d-1}_{\mathcal P}\ker(C_{0}\to M)=
\beta^{d-2}_{\mathcal P}\ker(C_{1}\to C_{0})=
\cdots = 
\beta^0_{\mathcal P}\ker(C_{d-1}\to C_{d-2}).
\]
Building a minimal resolution is an inductive procedure involving constructing minimal $\mathcal P$-covers of successive kernels, which typically requires a description of the entire $0^{\th}$ $\mathcal P$-Betti diagram of these kernels. 

One approach to calculating Betti diagrams that avoids this costly inductive procedure is the construction, given an object $M$ in $\mathcal{M}$ and an element $A$ in $\mathcal{P}$, of a chain complex $V = (\cdots \to V_1 \to V_0)$ of vector spaces such that, for all $d \geq 0$, $\dim H_d(V) = \beta^d_{\mathcal{P}} M(A)$.
Ideally, this construction is systematic, or even functorial in $M$.
This is the case for standard homological algebra for functors indexed by finite  upper semilattices, discussed in the next section, where such chain complexes are given by \emph{Koszul complexes}.
\end{step}

\section{Standard homological algebra}
\label{S:StdHomAlg}

Fix a field $K$, choose a finite poset $(J, \leqJ)$, and consider the abelian category $\Fun(J, \vect)$ whose objects are functors indexed by $J$ with values in the category $\vect$ of finite dimensional vector spaces. Morphisms in $\Fun(J, \vect)$ are the natural transformations. For an element $a$ in $J$, denote by $K(a, -) \colon J \to \vect$ the composition of the following functors, where the symbol $\text{set}$ denotes the category of finite sets:
\[
\begin{tikzcd}[column sep=50pt, row sep=2pt]
J \ar{r}{\hom_{J}(a, -)} & \text{set} \ar{r}{\text{free}} & \vect \\
& S \ar[mapsto, shorten=2ex]{r} & \bigoplus_{S} K
\end{tikzcd}
\]
That is, $K(a, b) = K$ if $b \geqJ a$ and $0$ otherwise, and $K(a, b \leqJ c) = \id_K$ if $b \geqJ a$ and $0$ otherwise. 
The functor $K(a, -)$ is called \emph{free on one generator in degree $a$}.

In this section we are going to describe the homological algebra of $\Fun(J, \vect)$, as presented in Section~\ref{S:RelHomAlg}, relative to the collection $\mathcal{S} := \{K(a, -) \mid a \in J\}$. This corresponds to standard homological algebra and all the following statements are well known: see for instance \cite{ARS1995, Weibel1994}.

\begin{step}[Independence]
The collection $\mathcal{S}$ is independent since its elements are indecomposable (see~\ref{Rel:Freeness}). This can be also shown 
using radicals, as in \cite{ARS1995}.
Recall that the \emph{radical} of a functor $F \colon J \to \vect$ is a subfunctor $\text{rad}(F) \subseteq F$ given by $\text{rad}(F)(a) = \sum_{s \ltJ a} \text{im}(F(s \ltJ a))$ for $a$ in $J$.
The quotient functor $H_0 F := F / \text{rad}(F)$ is \emph{semisimple}: that is, for every $a \ltJ b$ in $J$, the morphism $H_0 F(a \ltJ b)$ is the zero morphism. For $a$ and $b$ in $J$, $(H_0 K(a, -))(b)$ is $1$-dimensional if $a = b$ and $0$ otherwise.
Thus, if $F$ is free and isomorphic to $\bigoplus_{K(b,-) \in \mathcal{S}} K(b, -)^{\beta(K(b,-))}$, then $\beta(K(b,-)) = \text{dim}H_0 F(b)$ and hence the number $\beta(K(b,-))$ is uniquely determined by the isomorphism type of $F$.
\end{step}

\begin{step}[Exactness]
\label{Std:Exactness}
By the Yoneda lemma, for every functor $F \colon J \to \vect$, the linear map $\Nat_J(K(a, -),F) \to F(a)$, which assigns to a natural transformation $\varphi \colon K(a,-) \to F$ the element $\varphi_a(a \leqJ a)$ in $F(a)$, is a bijection.
Consequently, morphisms $f \colon F \to G$ and $g \colon G \to H$ in $\Fun(J,\vect)$ form an $\mathcal S$-exact sequence if, and only if,
for every $a$ in $J$, the linear maps $f_a\colon F(a)\to G(a)$ and $g_a\colon G(a)\to H(a)$ form an exact sequence of vector spaces. Thus, being $\mathcal{S}$-exact is the same as being exact. 
In particular, $f \colon F \to G$ is an $\mathcal S$-epimorphism if, and only if, it is a standard categorical epimorphism, which means $f_a \colon F(a) \to G(a)$ is surjective for all $a$ in $J$. To detect this we can again use the radical:
\end{step}

\begin{lemma}
\label{L:epi}
Let $J$ be a finite poset. A natural transformation $f \colon F \to G$ in $\Fun(J, \vect)$ is an epimorphism if, and only if, the composition of $f$ with the quotient $G \to G / \text{\rm rad}(G) = H_0 G$ is an epimorphism.
\end{lemma}
\begin{proof}
The only if part of this equivalence is clear. Conversely, suppose that the composition is an epimorphism and consider the set of all $a$ in $J$ for which $f_a$ is not a surjection. If this set is nonempty, then it contains a minimal element $a$ (since $J$ is finite).
Minimality of $a$ means that $\text{rad}(f)_a \colon \text{rad}(F)(a) \to \text{rad}(G)(a)$ is surjective.
This, combined with the surjectivity of the composition $F(a) \to G(a) \to H_0 G(a)$, implies that $f_a$ is surjective, which contradicts the assumption.
\end{proof}

Since $\mathcal S$-epimorphisms are the categorical epimorphisms, $\mathcal S$-projectives are the standard projectives, and it turns out that all of them are ($\mathcal{S}$-)free.

\begin{prop}
\label{L:StdAcyclic}
Let $J$ be a finite poset.
The collection $\mathcal{S}$ in $\Fun(J, \vect)$ is acyclic (see~\ref{Rel:Projectives}).
\end{prop}
\begin{proof}
We again use the radical.
Let $F$ be projective in $\Fun(J, \vect)$. Consider $H_0 F$ and the free functor $C_0 := \bigoplus_{a \in J} K(a, -)^{\dim H_0 F(a)}$. Note that $H_0 C_0$ and $H_0 F$ are isomorphic. We then use projectiveness to construct the following commutative diagram where all of the vertical arrows represent the quotient natural transformations:
\[
\begin{tikzcd}
F \ar{r} \ar{d} & C_0 \ar{d} \ar{r} & F \ar{d} \\
H_0F \ar[phantom, "\simeq"]{r} & H_0 C_0 \ar[phantom, "\simeq"]{r} & H_0F
\end{tikzcd}
\]
The vertical natural transformations in this diagram are epimorphisms, so by Lemma~\ref{L:epi} the horizontal natural transformations are as well.
Since the values of the functors are finite dimensional vector spaces, the horizontal natural transformations must therefore be isomorphisms.
\end{proof}

\begin{step}[\textbf{Minimal covers, resolutions, and Betti diagrams}]
\label{Std:MinCovers}
Let $F \colon J\to \vect$ be a functor. As in~\ref{L:StdAcyclic}, consider   $H_0 F = F / \text{rad}(F)$ and the free functor $C_0 := \bigoplus_{a \in J} K(a, -)^{\dim H_0 F(a)}$.
Since $H_0 C_0$ and $H_0 F$ are isomorphic, we can use the projectiveness
of $C_0$ to lift the natural transformation $C_0 \to H_0 F$ to a natural transformation $C_0 \to F$ which makes the following square commute, where the vertical arrows represent the quotient natural transformations: 
\[
\begin{tikzcd}
C_0 \ar{d} \ar{r} & F \ar{d} \\
H_0 C_0 \ar[phantom, "\simeq"]{r} & H_0 F
\end{tikzcd}
\]
The resulting natural transformation $C_0 \to F$ is an epimorphism (see discussion in~\ref{Std:Exactness}) and hence it is an $\mathcal S$-cover. It is also minimal by the same reasoning as in~\ref{L:StdAcyclic}.
It follows that there are enough $\mathcal S$-projectives and every functor $F \colon J \to \vect$ has a minimal $\mathcal S$-cover, and hence also a minimal $\mathcal S$-projective resolution.
We refer to $\mathcal S$-covers and $\mathcal S$-projective resolutions simply as \emph{covers} and \emph{resolutions}. 
Since every functor admits a minimal resolution, $\mathcal{S}$-Betti diagrams in all degrees are always defined for every functor. We also refer to them simply as \emph{Betti diagrams} and denote them as $\beta^d F \colon J \to \mathbb{N}$, omitting the symbol $\mathcal{S}$. 
A consequence of the above construction of a minimal cover of $F$ is the equality:
\[
\beta^0 F = \dim H_0 F.
\]
The sum $\sum_{a \in J} \beta^0 F(a)$ is called the \emph{number of generators of $F$}.

The Betti diagrams of $F$ are not independent from each other. For example, 
consider a differential in a minimal resolution of $F$:
\[\begin{tikzcd}
\overbrace{\displaystyle\bigoplus_{b\in
\supp(\beta^{d+1}F)}K(b,-)^{\beta^{d+1}F(b)}}^{C_{d+1}}\ar{r}{\partial} & 
\overbrace{\displaystyle\bigoplus_{a\in \supp(\beta^dF)}K(a,-)^{\beta^dF(a)}}^{C_d}
\end{tikzcd}\]
By minimality, each summand $K(b, -)$ of $C_{d+1}$
must map nontrivially onto at least one summand $K(a, -)$ of $C_d$. This means that every element in $\supp(\beta^{d+1}F)$
is bounded below by some element in $\supp(\beta^{d}F)$.
Relations between Betti diagrams are easier to describe when the indexing poset $J$ is an \emph{upper semilattice}.
In this case, to calculate Betti diagrams one can use Koszul complexes. Describing this is the content of the rest of this section.
\end{step}

\begin{step}[Upper semilattices]
\label{Std:UpperSemilattices}
Recall that the \emph{join} of a subset $S \subseteq J$ is an element in $J$, denoted by $\bigvee_J S$ or $\bigvee S$ when the poset is understood, satisfying the following universal property:
$s \leqJ \bigvee S$ for every $s$ in $S$, and, if $s \leqJ a$ in $J$ for every $s$ in $S$, then $\bigvee S \leqJ a$.
The join is the categorical coproduct in $J$.
Dually, the \emph{meet} of $S$ is an element in $J$, denoted by $\bigwedge_J S$ or $\bigwedge S$, satisfying the following universal property: $\bigwedge S \leqJ s$ for every $s$ in $S$, and, if $a \leqJ s$ in $J$ for every $s$ in $S$, then $a \leqJ \bigwedge S$. When they exist, joins and meets are unique.
If every nonempty subset of $J$ has a join, then $J$ is called an \emph{upper semilattice}.
If $J$ is an upper semilattice, then every subset $S$ which is bounded below (i.e., for which there exists $a$ in $J$ such that $a\leqJ s$ for every $s$ in $S$) also has a meet, given by $\bigvee \{a\in J\mid S \textnormal{ is bounded below by }  a\}$.

A subset $L \subseteq J$ of an upper semilattice is called a \emph{sublattice} if,
for every nonempty subset $S \subseteq L$, the join $\bigvee S$ belongs to $L$.
For $S \subseteq J$, we denote by $\langle S\rangle:=\left\{\bigvee T \mid T \neq \varnothing, T \subseteq S\right\} \subseteq J$ the smallest sublattice containing $S$.
\end{step}

\begin{step}[Koszul complexes]
\label{Std:KoszulComplexes}
Let $a$ be an element in $J$. Recall that an element $s$ in $J$ is called a \emph{parent} of $a$ if $s \ltJ a$ and there is no element $b$ in $J$ such that $s \ltJ b \ltJ a$. We also say that $a$ \emph{covers} $s$. We denote the set of parents of $a$ by $\mathcal{U}_J(a)$ or $\mathcal{U}(a)$ when the ambient poset is understood.

Suppose that $a$ has the following property: every subset $S \subseteq \mathcal{U}(a)$ that is bounded below 
admits a meet $\bigwedge S$. If $J$ is an upper semilattice, then every element $a \in J$ satisfies this property. 
Fix a total order $<$ on the set of parents $\mathcal{U}(a)$.
We then assign to every functor $F \colon J \to \vect$ a chain complex of $K$-vector spaces, called the \textbf{Koszul complex of $F$ at $a$}, defined for all $d \geq 0$ by
\[
(\mathcal{K}_a F)_d := \begin{cases}
F(a) & \text{ if } d = 0, \\
\displaystyle \bigoplus_{\substack{S \subseteq \mathcal{U}(a),\ |S| = d \\ S \text{ has lower bound}}} F({\textstyle\bigwedge S}) & \text{ if } d > 0.
\end{cases}
\]
For example, $(\mathcal{K}_a F)_1 = \bigoplus_{s \in \mathcal{U}(a)} F(s)$. We define the differentials as follows:
\begin{itemize}
\item For $d = 0$, define $\partial \colon (\mathcal{K}_a F)_1 \to (\mathcal{K}_a F)_0 = F(a)$ as the linear map which on the summand $F(s)$ in $\bigoplus_{s \in \mathcal{U}(a)} F(s) = (\mathcal{K}_a F)_1$ is given by $F(s \ltJ a)$.
\item For $d > 0$, define $\partial \colon (\mathcal{K}_a F)_{d + 1} \to (\mathcal{K}_a F)_d$ as the alternating sum $\partial = \sum_{i = 0}^d (-1)^i \partial_i$, where $\partial_i \colon (\mathcal{K}_a F)_{d + 1} \to (\mathcal{K}_a F)_d$ is the function mapping the summand $F(\bigwedge S)$ in $(\mathcal{K}_a F)_{d + 1}$, indexed by $S = \{s_0 < \cdots < s_{d}\} \subseteq \mathcal{U}(a)$, to the summand $F(\bigwedge (S \setminus \{s_i\}))$ in $(\mathcal{K}_a F)_{d}$, indexed by $S \setminus \{s_i\} \subseteq \mathcal{U}(a)$, via the function $F(\bigwedge S \leqJ \bigwedge (S \setminus \{s_i\}))$.
\end{itemize}
The linear functions $\partial$ form a chain complex as it is standard to verify that the composition of two consecutive such functions is null.

For a natural transformation 
$f \colon F \to G$, define:
\[
\begin{tikzcd}[column sep=40pt]
(\mathcal{K}_a F)_d \ar{r}{(\mathcal{K}_a f)_d} & (\mathcal{K}_a G)_d
\end{tikzcd}
:=
\begin{cases}
f_a & \text{ if } d = 0, \\
\displaystyle \bigoplus_{\substack{S \subseteq \mathcal{U}(a),\ |S| = d \\ S \text{ has lower bound}}} f_{\bigwedge S}
& \text{ if } d > 0.
\end{cases}
\]
These linear maps form a chain map denoted by $\mathcal{K}_a f \colon \mathcal{K}_a F \to \mathcal{K}_a G$.
The assignment $f \mapsto \mathcal{K}_a f$ is a functor from $\Fun(I, \vect)$ to the category of chain complexes called the \emph{Koszul complex at $a$}. Its fundamental property is the following.
\end{step}

\begin{thm}
\label{T:BettiKoszul}
Let $J$ be a finite poset. Let $a$ be an element in $J$ such that every nonempty subset $S \subseteq \mathcal{U}(a)$ 
that is bounded below
admits a meet $\bigwedge S$.
Then, for every functor $F \colon J \to \vect$ and every $d \geq 0$, the following equality holds:
\[
\beta^d F(a) = \dim H_d(\mathcal{K}_a F).
\]
\end{thm}

The assumption on the element $a$ in Theorem~\ref{T:BettiKoszul} is local, depending only on the parents of $a$. Under this assumption, which is satisfied for example if $J$ is an upper semilattice (see~\ref{Std:UpperSemilattices}), in order to calculate the Betti diagram $\beta^d F(a)$ of a functor $F$ at $a$, we do not need to construct the minimal resolution of $F$. Instead, it suffices to calculate the homology of the Koszul complex $\mathcal{K}_a F$, which can be done for each $a$  and each degree $d$ independently.

\begin{proof}[Proof of Theorem~\ref{T:BettiKoszul}]
The proof is based on three observations and follows closely what is presented in~\cite[Section~10]{tameness}. 
\smallskip

\noindent
\textbf{Observation 1.} 
For every natural transformation $f \colon F \to G$, the linear map $H_0(\mathcal{K}_a f)$ is isomorphic to $(H_0 f)_a$. 
In particular, the vector spaces $H_0(\mathcal{K}_a F)$ and $(H_0 F)(a)= (F / \text{rad}(F))(a)$ are isomorphic. Moreover, if $C_0 \to F$ is a minimal cover, then
$H_0(\mathcal{K}_a(C_0 \to F))$ is an isomorphism.
\smallskip

According to Observation 1, $\beta^0 F(a)= \dim (H_0 F)(a)$ (see~\ref{Std:MinCovers}) coincides with $\dim H_0(\mathcal{K}_a F)$, which is the statement of the theorem in the case $d = 0$. To prove this statement for $d>0$, we need two additional observations:
\smallskip

\noindent
\textbf{Observation 2.}
For $b$ in $J$:
\[
\dim H_d(\mathcal{K}_a K(b, -)) =
\begin{cases}
1 & \text{if } d = 0 \text{ and } a = b, \\
0 & \text{otherwise.}
\end{cases}
\]
Indeed, if $b \not\leqJ a$, then $\mathcal{K}_a K(b, -) = 0$ and the statement holds.
Otherwise, if $b = a$, then $(\mathcal{K}_a K(a, -))_0 = K$ and $(\mathcal{K}_a K(a, -))_d = 0$ for $d > 0$, and again the statement holds. 

Otherwise, $b \ltJ a$. Then $(\mathcal{K}_a K(b, -))_0 = K$ and, for every nonempty subset $S \subseteq \mathcal{U}(a)$ that is bounded below,
$K(b ,\bigwedge S)$ is $1$-dimensional if $S$ is a subset of $(b \leqJ \mathcal{U}(a)) := \{s \in \mathcal{U}(a) \mid b \leqJ s\}$, and $0$ otherwise. Consequently the complex $\mathcal{K}_a K(b,-)$
is isomorphic to
\[
\cdots \to \bigoplus_{\substack{S \subseteq (b \leqJ \mathcal{U}(a)) \\ |S| = 2}} K \to \bigoplus_{\substack{S \subseteq (b \leqJ \mathcal{U}(a)) \\ |S| = 1}} K \to K,
\]
which is the augmented chain complex of the standard $(|b \leqJ \mathcal{U}(a)| - 1)$-dimensional simplex, whose homology is trivial in all degrees. Thus in this case the statement also holds.
\smallskip

\noindent
\textbf{Observation 3.} 
Since taking direct sums is an exact operation, the Koszul complex at $a$ is an exact functor in the following sense: if $f \colon F \to G$ and $g \colon G \to H$ form an exact sequence, then the following is an exact sequence of chain complexes of vector spaces:
\[
\begin{tikzcd}
\mathcal{K}_a F
  \ar{r}{\mathcal{K}_a f}
&
\mathcal{K}_a G
  \ar{r}{\mathcal{K}_a g}
&
\mathcal{K}_a H.
\end{tikzcd}
\]

According to Observation 3, the Koszul complex functor at $a$ commutes with direct sums, which, combined with Observation 2, implies that for a free functor $C_0 = \bigoplus_{b \in J} K(b,-)^{\beta(b)}$,
\[
\dim H_d(\mathcal{K}_a C_0) =
\begin{cases}
\beta(a) & \text{if } d = 0, \\
0 & \text{otherwise.}
\end{cases}
\]

Let $F \colon J \to \vect$ be a functor and $C_0 \to F$ its minimal cover. This minimal cover fits into the following exact sequence:
\[
\begin{tikzcd}
0 \ar{r} & Z \ar{r} & C_0 \ar{r} & F \ar{r} & 0,
\end{tikzcd}
\]
which, according to Observation 3, leads to an exact sequence of chain complexes:
\[
\begin{tikzcd}
0 \ar{r} & \mathcal{K}_a Z \ar{r} & \mathcal{K}_a C_0 \ar{r} & \mathcal{K}_a F \ar{r} & 0,
\end{tikzcd}
\]
which in turn leads to a long exact sequence of vector spaces:
\[
\begin{tikzcd}
&&
\cdots
  \arrow[dll, phantom, ""{coordinate, name=Z}]
  \arrow[dll, rounded corners,
    to path={ -- ([xshift=2ex]\tikztostart.east)
    |- (Z)
    -| ([xshift=-2ex]\tikztotarget.west)
    -- (\tikztotarget)}]
\\
H_1(\mathcal{K}_a Z)
  \ar{r}
&
H_1(\mathcal{K}_a C_0)
  \ar{r}
&
H_{1}(\mathcal{K}_a F)
  \arrow[dll, phantom, ""{coordinate, name=Z}]
  \arrow[dll, rounded corners,
    to path={ -- ([xshift=2ex]\tikztostart.east)
    |- (Z)
    -| ([xshift=-2ex]\tikztotarget.west)
    -- (\tikztotarget)}]
\\
H_0(\mathcal{K}_a Z)
  \ar{r}
&
H_0(\mathcal{K}_a C_0)
  \ar{r}
&
H_{0}(\mathcal{K}_a F)
  \ar{r}
&
0.
\end{tikzcd}
\]
According to Observation 1, the function $H_0(\mathcal{K}_a C_0) \to H_{0}(\mathcal{K}_a F)$ is an isomorphism, and $H_0(\mathcal{K}_a Z)$ is isomorphic to $H_0(Z)$, whose dimension is $\beta^1 F(a)$.
Observation 2 gives the vanishing of $H_d(\mathcal{K}_a C_0)$ for $d \geq 1$.
Consequently, for $d\geq 1$, the map $H_{d}(\mathcal{K}_a F) \to H_{d - 1}(\mathcal{K}_a Z)$ is an isomorphism.
The case $d = 1$ then gives the equality between $\dim H_{0}(\mathcal{K}_a Z) = \beta^0 Z(a) = \beta^1 F(a)$ and $\dim H_{1}(\mathcal{K}_a F)$, which is the statement of the theorem for $d = 1$. The theorem for $d > 1$ follows by induction by applying what we already have proven to the functor $Z$.
\end{proof}

Theorem~\ref{T:BettiKoszul}, combined with the long exact sequence argument in its proof, implies the following.

\begin{cor}
\label{C:EqualBetti}
Let $J$ be a finite upper semilattice and $0\to F_0\to F_1\to F_2\to 0$ an exact sequence in $\Fun(J, \vect)$.
Let $a$ be an element in $J$ such that $\beta^d F_0(a) = 0$ for every $d \geq 0$. Then $\beta^d F_1(a) = \beta^d F_2(a)$ for every $d \geq 0$.
\end{cor}

\section{Filtrations, subfunctors, and their standard Betti diagrams}
\label{S:FiltrationsSubfunctors}

In this section, we discuss in more depth resolutions and Betti diagrams of 
certain filtrations. Recall that a functor $F \colon J \to \vect$ is called a \emph{filtration} if, for all $a \leqJ b$ in $J$, $F(a \leqJ b)$ is a monomorphism.
For example, free functors and constant functors are filtrations, and in particular the constant functor $K_J \colon J \to \vect$, whose values are the $1$-dimensional vector space $K$, is a filtration.
Since subfunctors of a filtration are filtrations, the kernel of a cover $C_0 \to F$ is a filtration.
Two results of this section are of particular interest to us.
One is Corollary~\ref{C:FiltrationBettiContainments}, which describes relations between the supports of Betti diagrams of a filtration whose $0^\th$ Betti diagram has support bounded below.
The other one is Corollary~\ref{C:SubfunctorBettiContainments}, which describes a similar, but weaker, relation for any subfunctor of the constant functor $K_J$, even if its $0^{\text{th}}$ Betti diagram does not have a support bounded below.

The following theorem follows from results in~\cite{tameness}
(see also~\cite[Lemma 2.1]{Vipond2020}).

\begin{thm}
\label{T:SublatticeDiscretization}
Let $J$ be a finite upper semilattice and $F\colon J\to\text{\rm vect}_K$ a functor.
 Then the following containment holds:
\[
\bigcup_{d\geq 0} \supp(\beta^{d}F) \subseteq \langle \supp(\beta^{0}F) \cup \supp(\beta^{1}F) \rangle.
\]
\end{thm}
\begin{proof}
Consider a minimal presentation of $F$ given by an exact sequence:
\[
\begin{tikzcd}
\displaystyle \bigoplus_{a \in \supp(\beta^1F)} \!\!\!\! K(a, -)^{\beta^1 F(a)}
 \ar{r}{\partial}
& 
\displaystyle \bigoplus_{a \in \supp(\beta^0F)} \!\!\!\! K(a, -)^{\beta^0 F(a)}
 \ar{r}
&
F
 \ar{r}
&
0,
\end{tikzcd}
\]

Using the vocabulary of \cite{tameness}, we can discretize the functors of this sequence by the sublattice $\langle \supp(\beta^{0}F)\cup \supp(\beta^{1}F)\rangle$. By \cite[Corollary~10.18.(2)]{tameness}, we conclude that the Betti diagrams of $F$ indexed by $J$ are isomorphic to those of $F$ restricted to the sublattice
$\langle \supp(\beta^{0}F)\cup \supp(\beta^{1}F)\rangle\subset J$, and consequently
\[
\bigcup_{d\geq 0} \supp(\beta^{d}F)\subseteq \langle \supp(\beta^{0}F)\cup \supp(\beta^{1}F)\rangle.
\qedhere
\]
\end{proof}

As a consequence, if $\beta^dF(a)\neq0$, then there is a subset $S \subseteq \supp(\beta^{0}F) \cup \supp(\beta^{1}F)$ for which $a=\bigvee S$. This gives a considerable restriction on what elements of the upper semilattice $J$ can belong to $\supp(\beta^{d}F)$ for $d>1$.

\begin{cor}
\label{C:BettiContainments}
Let $J$ be a finite upper semilattice and $F \colon J \to \vect$ a functor for which $\supp(\beta^{0} F)$ is bounded below (the meet $\bigwedge \supp(\beta^{0} F)$ exists). Then
\[
\cdots
\subseteq \langle \supp(\beta^{3} F) \rangle
\subseteq \langle \supp(\beta^{2} F) \rangle
\subseteq \langle \supp(\beta^{1} F) \rangle.
\]
\end{cor}
\begin{proof}
Let $b_0 = \bigwedge \supp(\beta^{0}F)$.
For all $a$ in $\supp(\beta^{0}F)$, we have $a \geqJ b_0$, so there is a monomorphism $K(a,-) \subseteq K(b_0,-)$. These monomorphisms, for all $a$ in $\supp(\beta^{0}F)$, fit into the following pushout square where the top horizontal arrow represents a minimal cover of $F$:
\[
\begin{tikzcd}
\displaystyle \bigoplus_{a \in \supp(\beta^{0} F)} \!\!\!\! K(a, -)^{\beta^0 F(a)}
 \ar{r}
 \ar[hook]{d}
&
F
 \ar[hook]{d}
\\
\displaystyle \bigoplus_{a \in \supp(\beta^{0} F)} \!\!\!\! K(b_0, -)^{\beta^0 F(a)}
 \ar{r}
&
G
\end{tikzcd}
\]
Then the natural transformation represented by the bottom horizontal arrow is a minimal cover of $G$.
Furthermore, since the kernels of the top and bottom horizontal natural transformations coincide, we get:
\[
\beta^d G(a) = \begin{cases}
\sum_{b \in J} \beta^0 F(b) & \text{if } d = 0 \text{ and } a = b_0, \\
0 & \text{if } d = 0 \text{ and } a \neq b_0, \\
\beta^d F(a) & \text{if } d > 0.
\end{cases}
\]
In particular $\supp(\beta^{0} G) = \{b_0\}$ and $\supp(\beta^{d} G) = \supp(\beta^{d} F)$ for every $d \geq 1$.
Since every element in $\supp(\beta^{2} G)$ is bounded by some element in $\supp(\beta^{1} G)$, and every element in $\supp(\beta^{1} G)$ is bounded by some element in $\supp(\beta^{0} G) = \{b_0\}$, we can use Theorem~\ref{T:SublatticeDiscretization} to conclude that $\supp(\beta^{2} G)\subseteq \langle \supp(\beta^{1} G) \rangle$ and consequently $\langle \supp(\beta^{2} G) \rangle \subseteq \langle \supp(\beta^{1} G) \rangle$. We then proceed by induction on $d$.
\end{proof}

For a filtration, Corollary~\ref{C:BettiContainments} can be extended to include also the $0^\th$ Betti diagram.

\begin{cor}
\label{C:FiltrationBettiContainments}
Let $J$ be a finite upper semilattice and $F \colon J \to \vect$ a filtration such that $\supp(\beta^{0}F)$ is bounded below. Then
\[
\cdots
\subseteq \langle \supp(\beta^{2} F) \rangle
\subseteq \langle \supp(\beta^{1} F) \rangle
\subseteq \langle \supp(\beta^{0} F) \rangle.
\]
\end{cor}
\begin{proof}
The assumptions on $F$ are equivalent to $F$ being a subfunctor of $G := K(b, -)^n$ for some $b$ in $J$ and a natural number $n$.
Consider the exact sequence $F \hookrightarrow G \to G / F$. Since $G$ is free, $\beta^1(G / F)(a) \leq \beta^0 F(a)$ and $\beta^{d + 1}(G / F)(a) = \beta^d F(a)$ for $d > 0$ and  every  $a$ in $J$.
Consequently, $\supp(\beta^1(G / F)) \subseteq \supp(\beta^0 F)$ and $\supp(\beta^{d + 1}(G / F)) = \supp(\beta^d F)$ for $d > 0$. By applying Corollary~\ref{C:BettiContainments} to the quotient $G / F$, we obtain the desired
containments. For example,
\[
\langle \supp(\beta^1 F) \rangle = \langle \supp(\beta^{2}(G / F)) \rangle
\subseteq \langle \supp(\beta^{1}(G / F)) \rangle
\subseteq \langle \supp(\beta^0 F) \rangle.
\qedhere
\]
\end{proof}
\begin{step}[Subfunctors]
Next we focus on subfunctors of the constant functor $K_J \colon J \to \vect$.  Given a subfunctor $F$ of $K_J$, we denote by $F \subseteq K_J$ the inclusion natural transformation.
Note that a functor $F \colon J \to \vect$ is isomorphic to a subfunctor of $K_J$ if, and only if, $F$ is a filtration and $\dim F(a) \leq 1$ for every $a$ in $J$ (compare with assumptions in Corollary~\ref{C:FiltrationBettiContainments}). We denote by $\text{Sub}(J)$ the set of all  subfunctors of $K_J$.
We consider two ways of parameterizing  
this set by posets: via upsets and antichains.
\end{step}

\begin{step}[Upsets]
\label{D:SubfunctorsUpsets}
An \emph{upset} of a poset $(J,\leqJ)$ is a subset $U \subseteq J$ such that, for all $a$ in $U$ and $b \geqJ a$, the element $b$ is in $U$. The set of upsets of $J$ is denoted $\text{Up}(J)$ and comes naturally equipped with a distributive lattice structure where the order relation is the inclusion $\subseteq$, joins are unions, and meets are intersections. An example of an upset is the \emph{support} of a subfunctor $F \subseteq K_J$, defined as the subset $\supp(F) := \{a \in J \mid F(a) \neq 0\}$. In fact, the function $F \mapsto \supp(F)$ is a bijection between the collections of subfunctors of $K_J$ and upsets in $J$. Its inverse sends an upset $U$ to the subfunctor $K_U \subseteq K_J$ defined as, for all $a$ in $J$,
\begin{align*}
K_U(a) & = \begin{cases}
K  & \text{if } a \in U, \\
0 & \text{otherwise.}
\end{cases}
\end{align*}
Thus $K_U \subseteq K_J$ is the unique subfunctor for which 
$\text{supp}(K_U)=U$.
The inclusion of subfunctors $F \subseteq G \subseteq K_J$ defines a poset relation that makes this bijection a poset isomorphism, between $(\text{Sub}(J),\subseteq)$ and 
$(\text{Up}(J),\subseteq)$. We are however more interested in the opposite poset to  $(\text{Sub}(J),\subseteq)$. This is because of its relation to antichains.
\end{step}

\begin{step}[Antichains]
\label{D:Antichains}
An \emph{antichain} of a poset $(J,\leqJ)$ is a subset of $J$ whose elements are  pairwise incomparable. In particular, singletons of $J$ are antichains. The set of antichains of $J$ is denoted by $\text{Anti}(J)$. Given an upset $U$ of $J$, the set $\Min(U)$ of its minimal elements is an antichain. In fact, the function $U \mapsto \Min(U)$ from upsets to antichains is a bijection, and so we could equip antichains with the induced poset structure from $\text{Up}(J)$. 
However, we would like for the subposet of singletons in 
$\text{Anti}(J)$ to coincide with the poset $J$. For this reason, we define the relation on antichains by the opposite order on upsets. Thus, given two upsets $U$ and $V$, we write $\Min(U) \leqJ \Min(V)$ if $V \subseteq U$.
Explicitly, for two antichains $S$ and $T$ in $J$,
the relation  $S \leqJ T$ holds if, and only if, for every
$t$ in $T$, there is $s$ in $S$ such that $s \leqJ t$ in $J$.
In this way, $(\text{Anti}(J), \leqJ)$ is a distributive lattice whose restriction to singletons coincides with $J$.

To a subfunctor $F \subseteq K_J$ we associate the antichain $\Min \supp(F)$, which coincides with $\supp(\beta^0 F)$. The function $F \mapsto \supp(\beta^0 F)$ is a bijection
between the collection $\text{Sub}(J)$ of subfunctors of $K_J$ and the set of antichains $\text{Anti}(J)$. Its inverse sends an antichain $S$ to the subfunctor $K(S, -) := K_{(S \leqJ J)} \subseteq K_J$ where $(S \leqJ J)$ denotes the upset $\{a \in J \mid \exists s \in S, s \leqJ a\}$. This induces a distributive lattice structure on subfunctors of $K_J$: given two subfunctors $F, G \subseteq K_J$, we write $F \leqJ G$ if $\supp(\beta^0 F) \leqJ \supp(\beta^0 G)$, which is equivalent to the inclusion $\supp(F) \supseteq \supp(G)$. We denote this poset by $(\text{Sub}(J), \leqJ)$.
\end{step}

\begin{step}[Global Koszul complexes]
We now discuss free resolutions of subfunctors of $K_J$ under the assumption that $J$ is an upper semilattice. Let $F \subseteq K_J$ be a subfunctor. Fix a total order $<$ on the antichain $\supp(\beta^0 F)$.
This order is used to construct a chain complex called \textbf{global Koszul complex} of $F$. For $d \geq 0$, define:
\[
(\mathcal{K}F)_d := 
\bigoplus_{\substack{S \subseteq \supp(\beta^0 F) \\ |S| = d + 1}} K({\textstyle\bigvee S, -}).
\]
For example, $(\mathcal{K} F)_0 = \bigoplus_{s \in \supp(\beta^0 F)} K(s, -)$, which coincides with the minimal cover of $F$ (see~\ref{Std:MinCovers}).
Set $\partial \colon (\mathcal{K} F)_{d + 1} \to (\mathcal{K} F)_d$ to be the alternating sum $\partial = \sum_{i = 0}^{d + 1} (-1)^i \partial_i$, where $\partial_i \colon (\mathcal{K} F)_{d + 1} \to (\mathcal{K} F)_d$ is the function mapping the summand $K(\bigvee S, -)$ in $(\mathcal{K} F)_{d + 1}$, indexed by $S = \{s_0 < \cdots < s_{d + 1}\} \subseteq \supp(\beta^0F)$, to the summand $K(\bigvee (S \setminus \{s_i\}), -)$ in $(\mathcal{K} F)_{d}$, indexed by $S \setminus \{s_i\} \subseteq \supp(\beta^0 F)$, via the inclusion $K(\bigvee S, -) \subseteq K(\bigvee (S \setminus \{s_i\}),-)$. 
Finally define $\partial \colon (\mathcal{K} F)_0\to F$ to be a minimal cover of $F$. It is standard to verify that the natural transformations $\partial$ form a chain complex.
\end{step}

\begin{prop}
Let $J$ be a finite upper semilattice and $F \subseteq K_J$ be a subfunctor.
Then $\mathcal{K} F \to F$ is a free resolution of $F$.
\end{prop}
\begin{proof}
We need to show that, for every $a$ in $J$, the complex of vector spaces
$(\mathcal{K} F)(a) \to F(a)$ is exact. Consider the set $T: = \{s \in \supp(\beta^0 F) \mid s \leqJ a\}$.
If this set is empty, then the complex $(\mathcal{K} F)(a) \to F(a)$ is trivial and hence it is exact. 
Otherwise, the complex $(\mathcal{K} F)(a) \to F(a)$
is isomorphic to
\[
\cdots \to \bigoplus_{\substack{S \subseteq T \\ |S| = 2}} K \to \bigoplus_{\substack{S \subseteq T \\ |S| = 1}} K \to K,
\]
which is the augmented chain complex of the standard $(|T| - 1)$-dimensional simplex, whose homology is trivial in all degrees.
\end{proof}

The global Koszul complex $\mathcal{K} F \to F$ is a free resolution of $F$ which may fail however to be minimal. A minimal resolution of $F$ is a direct summand of the global Koszul complex, which gives:

\begin{cor}
\label{C:SubfunctorBettiContainments}
Let $J$ be a finite upper semilattice and $F \subseteq K_J$ be a subfunctor.
Then the following containment holds:
\[
\bigcup_{d\ge 0}\supp(\beta^d F) \subseteq \langle \supp(\beta^0 F) \rangle.
\]
\end{cor}

\section{Constructing relative homological algebra}
\label{S:GradedRelHomAlg}

Let $(I, \leq)$ be a finite poset. In this section we present a strategy for defining an independent (see~\ref{Rel:Freeness}) and acyclic (see~\ref{Rel:Projectives}) collection for relative homological algebra in the category $\Fun(I,\vect)$.
Our starting point and the standard assumption 
throughout the entire section is a functor $\parameterization \colon J^{\op} \to \Fun(I, \vect)$ where  $(J,\leqJ)$ is a finite poset.

We define the collection $\mathcal{P} := \{\parameterization(a) \mid a \in J,\ \parameterization(a) \neq 0\}$. Although no structure on the collection $\mathcal{P}$ was needed in our description of the homological algebra relative to $\mathcal{P}$ in Section~\ref{S:RelHomAlg}, in this section we illustrate the advantage of having the grading on $\mathcal{P}$ given by the poset structure on $J$.
The functor $\parameterization$ is used to translate between the $\mathcal{P}$-relative homological algebra on $\Fun(I, \vect)$ and the standard homological algebra on $\Fun(J, \vect)$.
This translation is done via a pair of adjoint functors.

\begin{example}
\label{E:Running}
Throughout this section, we  consider the following running example. Let $(I, \leq)$ be the set $\{0, 1\}^2$ with the product order, $J$ the set $\{(u, v) \in I^2 \mid u \leq v\}$ with the product order, and
\[
\parameterization \colon \left\{ \begin{array}{ccc}
J^{\op} & \to & \Fun(I, \vect) \\
(u, v) & \mapsto & \coker(K(v, -) \to K(u, -))
\end{array}
\right.
\]
the translation or \emph{parameterization} functor.
Explicitly, we have the Hasse diagrams (where the order goes up and to the right)
\[
I = \begin{tikzpicture}[every node/.append style={font=\scriptsize, inner sep=1pt}, baseline = {([yshift = -.5ex, scale=1](0, .5))}]
\node (00) at (0, 0) {$00$};
\node (01) at (0, 1) {$01$};
\node (10) at (1, 0) {$10$};
\node (11) at (1, 1) {$11$};

\draw[dotted] (00) -- (01) -- (11) -- (10) -- (00) -- cycle;
\end{tikzpicture}
\qquad \textnormal{and} \qquad
J =
\begin{tikzpicture}[every node/.append style={font=\scriptsize, inner sep=1pt}, baseline = {([yshift = -.5ex, scale=1](0, 1))}]
\node (0000) at (0, 0) {$00, 00$};
\node (0010) at (1, 0) {$00, 10$};
\node (1010) at (2, 0) {$10, 10$};
\node (0001) at (0, 1) {$00, 01$};
\node (0101) at (0, 2) {$01, 01$};
\node (0011) at (1, 1) {$00, 11$};
\node (1011) at (2, 1) {$10, 11$};
\node (0111) at (1, 2) {$01, 11$};
\node (1111) at (2, 2) {$11, 11$};

\draw[dotted] (0001) -- (0000) -- (0010) -- (1010) -- (1011) -- (0011) -- (0001) -- (0101) -- (0111) -- (1111) -- (1011) -- (1010)
(0111) -- (0011) -- (0010);
\end{tikzpicture}.
\]
The image of $J$ by $\parameterization$ is then
\[
\begin{array}{ccc}
0
&
\begin{tikzpicture}[scale = \diagramscale, vcenter=1]
\tikzgrid{1}
\fill[functor fill, functor draw] (0, 1) rectangle (1, 2);
\useasboundingbox (0, 0) rectangle (2, 2);
\end{tikzpicture}
&
0
\\[2ex]
\begin{tikzpicture}[scale = \diagramscale, vcenter=1]
\tikzgrid{1}
\fill[functor fill, functor draw] (0, 0) rectangle (2, 1);
\end{tikzpicture}
&
\begin{tikzpicture}[scale = \diagramscale, vcenter=1]
\tikzgrid{1}
\fill[functor fill, functor draw] (0, 0) -- (2, 0) -- (2, 1) -- (1, 1) -- (1, 2) -- (0, 2) -- cycle;
\useasboundingbox (0, 0) rectangle (2, 2);
\end{tikzpicture}
&
\begin{tikzpicture}[scale = \diagramscale, vcenter=1]
\tikzgrid{1}
\fill[functor fill, functor draw] (1, 0) rectangle (2, 1);
\end{tikzpicture}
\\[2ex]
0
&
\begin{tikzpicture}[scale = \diagramscale, vcenter=1]
\tikzgrid{1}
\fill[functor fill, functor draw] (0, 0) rectangle (1, 2);
\useasboundingbox (0, 0) rectangle (2, 2);
\end{tikzpicture}
&
0
\end{array},
\]
where we have represented each nonzero functor $\parameterization(a) \colon I \to \vect$ in $\mathcal{P}$ by their support: the dotted squares correspond to the poset $I$, the shaded vertices correspond to the support of the functor, where the functor is equal to $K$, and the transition functions are all implicitly of maximal rank. The position of each functor $\parameterization(a)$ corresponds to the position of $a$ in the Hasse diagram of $J$. 

Finally, we have the collection
\[
\mathcal{P} = \left\{
\begin{tikzpicture}[scale = \diagramscale, vcenter=1]
\tikzgrid{1}
\fill[functor fill, functor draw] (0, 0) rectangle (1, 2);
\end{tikzpicture}
\; , \;
\begin{tikzpicture}[scale = \diagramscale, vcenter=1]
\tikzgrid{1}
\fill[functor fill, functor draw] (0, 0) rectangle (2, 1);
\end{tikzpicture}
\; , \;
\begin{tikzpicture}[scale = \diagramscale, vcenter=1]
\tikzgrid{1}
\fill[functor fill, functor draw] (0, 0) -- (2, 0) -- (2, 1) -- (1, 1) -- (1, 2) -- (0, 2) -- cycle;
\end{tikzpicture}
\; , \;
\begin{tikzpicture}[scale = \diagramscale, vcenter=1]
\tikzgrid{1}
\fill[functor fill, functor draw] (1, 0) rectangle (2, 1);
\end{tikzpicture}
\; , \;
\begin{tikzpicture}[scale = \diagramscale, vcenter=1]
\tikzgrid{1}
\fill[functor fill, functor draw] (0, 1) rectangle (1, 2);
\end{tikzpicture}
\;
\right\}.
\]
This example is a specific instance of \emph{lower hooks}, introduced in \cite{BOO2021}, which we study in general in \ref{E:LowerHooks}.
\end{example}

\begin{step}[Adjunction]
\label{D:Adjunction}
The functor $\parameterization$ induces the following pair of adjoint functors:
\[
\begin{tikzcd}
\Fun(I,\vect)
 \ar[bend left=10]{r}{\mathcal R}
&
\Fun(J,\vect)
 \ar[bend left=10]{l}{\mathcal L}
\end{tikzcd}
\]
where
\begin{itemize}[leftmargin=2em]
\item $\mathcal R := \Nat_I(\parameterization, -)$ assigns to $M \colon I \to \vect$ the functor $\Nat_I(\parameterization(-), M) : J \to \vect$,
\item $\mathcal L$ assigns to $F\colon J\to\vect$ the following colimit in $\Fun(I,\vect)$:
\[
\mathcal L F := \text{colim} \left( \begin{tikzcd}
\displaystyle \bigoplus_{a_0 \ltJ a_1 \text{ in } J} \parameterization(a_1) \otimes F(a_0)
 \ar[shift left=5pt]{r}{d_0}
 \ar[shift right=5pt]{r}[swap]{d_1}
&
\displaystyle \bigoplus_{a\text{ in } J} \parameterization(a) \otimes F(a)
\end{tikzcd} \right),
\]
where the summand $\parameterization(a_1) \otimes F(a_0)$, indexed by $a_0 \ltJ a_1$, is mapped
\begin{itemize}[leftmargin=1em]
\item by $d_0$ to the summand $\parameterization(a_1) \otimes F(a_1)$ via the morphism $\text{id}_{\parameterization(a_1)} \otimes F(a_0 \ltJ a_1)$;
\item by $d_1$ to the summand $\parameterization(a_0) \otimes F(a_0)$ via the morphism $\parameterization(a_0 \ltJ a_1) \otimes \text{id}_{F(a_0)}$.
\end{itemize}
\end{itemize}
\end{step}

\begin{prop}
\label{sadgsghj}
The functor $\mathcal L$ is left adjoint to $\mathcal R$.
\end{prop}
\begin{proof}
Let $F \colon J \to \vect$ and $M \colon I \to \vect$ be two functors. By the universal property of the colimit and the tensor-hom adjunction, the set $\Nat_I(\mathcal{L}F, M)$ is in natural bijection with the set of sequences of linear maps $\{f_a \colon F(a) \to \Nat_I(\parameterization(a), M)\}_{a\in J}$ making the following square commute for every $a_0 \ltJ a_1$ in $J$:
\[
\begin{tikzcd}
F(a_0)
 \ar{r}{f_{a_{0}}}
 \ar{d}[swap]{F(a_0 \ltJ a_1)}
&
\Nat_I(\parameterization(a_0), M)
 \ar{d}{\Nat_I(\parameterization(a_0 \ltJ a_1), M)}
\\
F(a_1)
 \ar{r}{f_{a_1}}
&
 \Nat_I(\parameterization(a_1), M)
\end{tikzcd}
\]
Thus this set of sequences is in bijection with $\Nat_J(F, \Nat_I(\parameterization, M)) = \Nat_J(F, \mathcal{R} M)$.
In this way we get a natural isomorphism between $\Nat_I(\mathcal{L} F, M)$ and $\Nat_J(F, \mathcal{R} M)$, which gives the desired adjunction. 
\end{proof}

Let $a$ be in $J$ and consider the free functor $K(a,-) \colon J \to \vect$. By the adjunction between $\mathcal L$ and $\mathcal R$, for every $M \colon I \to \vect$, the set $\Nat_I(\mathcal{L} K(a,-), M)$ is naturally isomorphic to the set $\Nat_J(K(a, -), \mathcal{R} M)$, which is in natural bijection with $\mathcal{R} M(a) = \Nat_I(\parameterization(a), M)$.
Consequently $\mathcal{L} K(a, -)$ and $\parameterization(a)$ are isomorphic functors for every $a$ in $J$, and the collection $\mathcal P=\{\parameterization(a) \mid a \in J,\ \parameterization(a) \neq 0\}$ can be identified with $\{\mathcal{L} K(a, -) \mid a \in J,\ \mathcal{L}K(a,-) \neq 0\}$.

\begin{example}
\label{E:M}
Continuing  our running Example~\ref{E:Running}, consider the functor $M \colon I \to \vect$ represented as follows:
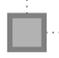
\[
\begin{tikzpicture}[scale = \diagramscale, vcenter=1]
\tikzgrid{1}
\fill[functor fill, functor draw] (0, 0) rectangle (1, 1);
\end{tikzpicture}
\]
Then, using the same representation conventions for functors indexed by $J$, the functor $\mathcal{R} M$ is represented as
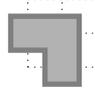
\[
\begin{tikzpicture}[scale = \diagramscale, vcenter=2]
\tikzgrid{2}
\fill[functor fill, functor draw] (0, 1) -- (1, 1) -- (1, 0) -- (2, 0) -- (2, 2) -- (0, 2) -- cycle;
\end{tikzpicture}
\]
\end{example}

\begin{step}[Adjunction unit and counit]
For $M \colon I \to \vect$, the \emph{counit} $\varepsilon_M \colon \mathcal{LR} M \to M$ is the natural transformation adjoint to $\text{id} \colon \mathcal{R} M \to \mathcal{R} M$.
For $F \colon J \to \vect$, the \emph{unit} $\eta_F\colon F\to \mathcal{RL} F$ is the natural transformation adjoint to $\text{id} \colon \mathcal{L} F \to \mathcal{L} F$. If $F = K(a,-)$ for some $a$ in $J$, then $\eta_F$ is also denoted by $\eta_a$. 
The adjunction between $\mathcal{L}$ and $\mathcal{R}$ implies the commutativity of the following diagrams, which in particular implies that $\mathcal{R} \varepsilon_{M}$ and $\varepsilon_{\mathcal{L}F}$ are epimorphisms, and $\eta_{\mathcal{R} M} $ and $\mathcal{L} \eta_F$ are monomorphisms: 
\[
\begin{tikzcd}
\mathcal{R} M
 \ar{r}{\eta_{\mathcal{R} M}}
 \ar[bend right, swap]{dr}{\text{id}}
&
\mathcal{RLR} M
 \ar{d}{\mathcal{R} \varepsilon_{M}}
\\
&
\mathcal{R} M
\end{tikzcd}
\hspace{10mm}
\begin{tikzcd}
\mathcal{LRL} F
 \ar[swap]{d}{\varepsilon_{\mathcal{L} F}}
&
\mathcal{L}F \ar{l}[swap]{\mathcal{L} \eta_F}
 \ar[bend left, swap]{dl}[swap]{\text{id}}
\\
\mathcal{L} F
\end{tikzcd}
\]
\end{step}

The functor $\mathcal{R}$ can be used to translate
between $\mathcal{P}$-exactness in $\Fun(I, \vect)$ and the standard exactness in $\Fun(J, \vect)$:

\begin{prop}
\label{P:RelAndStdExactness}
Let $I$ and $J$ be finite posets.
A pair of natural transformations $f \colon M \to N$ and $g \colon N \to L$ in $\textnormal{Fun}(I, \vect)$ form a $\mathcal{P}$-exact sequence (see~\ref{Rel:Exactness}) if, and only if, the natural transformations $\mathcal{R} f$ and $\mathcal{R} g$ form an exact sequence in $\textnormal{Fun}(J, \textnormal{vect}_K)$.
\end{prop}
\begin{proof}
By definition of the functor $\mathcal{R} = \Nat_{I}(\parameterization, -)$, the pair $\mathcal{R} f$ and $\mathcal{R} g$ forms an exact sequence if, and only if, $\Nat_{I}(\parameterization(a), f)$ and $\Nat_{I}(\parameterization(a), g)$ form an exact sequence for all $a$ in $J$, which corresponds to $\mathcal{P}$-exactness in $\Fun(I, \vect)$.
\end{proof}

\begin{cor}
\label{C:RelAndStdEpi}
A natural transformation $f \colon M \to N$ in $\Fun(I,\vect)$ is a $\mathcal P$-epimorphism if, and only if, $\mathcal{R} f$ is an epimorphism in $\Fun(J, \vect)$.
\end{cor}

For example, consider a functor $M \colon I \to \vect$ and the natural transformation $\varepsilon_M \colon \mathcal{LR} M \to M$. Since $\mathcal{R} \varepsilon_M$ is an epimorphism, $\varepsilon_M$ is therefore a $\mathcal{P}$-epimorphism. 

To translate between minimal $\mathcal P$-covers in $\Fun(I, \vect)$ and standard minimal covers in $\Fun(J, \vect)$, an additional assumption on the functor $\parameterization \colon J^{\op} \to \Fun(I, \vect)$ is required.
Here is the key definition of our paper.

\begin{step}[Thinness]
\label{D:Thinness}
The functor $\parameterization \colon J^{\op} \to \Fun(I, \vect)$ is called \emph{thin} if, for every $a$ in $J$, the unit natural transformation $\eta_a \colon K(a, -) \to \mathcal{RL} K(a,-)$ is an epimorphism.
\end{step}

Our primary examples of thin functors are given by the following result:

\begin{prop}
\label{P:SingleSourceThinness}
Let $I$ and $J$ be finite posets and $\parameterization \colon J^\op \to \Fun(I, \vect)$ a functor. 
Suppose that
\begin{itemize}
\item for all $a$ in $J$, the functor $\parameterization(a)$ has at most one generator, i.e.\ $\sum_{v\in I} \beta^0 \parameterization(a)(v) \le 1$.
\item for all $a, b$ in $J$, $\Nat(\parameterization(b), \parameterization(a)) \neq 0$ only if $a \leqJ b$. 
\end{itemize}
Then the functor $\parameterization$ is thin.
\end{prop}

\begin{proof}
Let $a$ be in $J$. If $\parameterization(a) = 0$, then $\mathcal{RL} K(a, -) = \Nat(\parameterization(-), \parameterization(a)) = 0$, so the unit natural transformation $\eta_a \colon K(a, -) \to \mathcal{RL} K(a, -)$ is surjective.

Otherwise, let $b$ be another element of $J$. If $\parameterization(b) = 0$, then $\Nat(\parameterization(b), \parameterization(a)) = 0$. If $a \not\leqJ b$, then by hypothesis we also have $\Nat(\parameterization(b), \parameterization(a)) = 0$. Otherwise, we have $a \leqJ b$ and $\parameterization(b) \neq 0$. Using the hypothesis on $\parameterization$, write $\supp(\beta^0 \parameterization(b)) = \{x_b\}$. Every natural transformation $\varphi \colon \parameterization(b) \to \parameterization(a)$ is then determined by $\varphi_{x_b} \colon \parameterization(b)(x_b) \to \parameterization(a)(x_b)$. Since $\dim\parameterization(b)(x_b) = 1$ and $\dim\parameterization(a)(x_b)\leq 1$, we get $\dim \Nat(\parameterization(b), \parameterization(a)) \leq 1$.
Moreover, since $a \leqJ b$, we have $\dim K(a, b) = 1$, and so the natural map $\eta_a(b) \colon K(a, b) \to \Nat(\parameterization(b), \parameterization(a))$ is a surjection.
Thus the unit natural transformation $\eta_a \colon K(a, -) \to \mathcal{RL} K(a, -) = \Nat(\parameterization(-), \parameterization(a))$ is surjective, and we conclude that $\parameterization$ is thin.
\end{proof}

\begin{example}
\label{E:ManualThinness}
In our running Example~\ref{E:Running}, the functor $\parameterization$ is thin. Lower hooks fit into the setting of Proposition \ref{P:SingleSourceThinness}, but we can also check thinness by hand: for instance, for $a = (00, 10)$ in $J$, we have
\[
\begin{tikzpicture}[scale = \diagramscale, vcenter=2]
\tikzgrid{2}
\fill[functor fill, functor draw] (1, 0) rectangle (3, 3);
\end{tikzpicture} = K(a, -)
\; \xrightarrow{\eta_a} \;
\mathcal{RL} K(a, -) = \mathcal{R}\,
\begin{tikzpicture}[scale = \diagramscale, vcenter=1]
\tikzgrid{1}
\fill[functor fill, functor draw] (0, 0) rectangle (1, 2);
\end{tikzpicture}
=
\begin{tikzpicture}[scale = \diagramscale, vcenter=2]
\tikzgrid{2}
\fill[functor fill, functor draw] (1, 0) rectangle (2, 3);
\end{tikzpicture}\;,
\]
so the unit natural transformation $\eta_a \colon K(a, -) \to \mathcal{RL} K(a, -)$ is indeed an epimorphism.
\end{example}

Among thin functors there are functors that satisfy a stronger requirement.

\begin{step}[Flatness]
\label{D:Flatness}
The functor $\parameterization \colon J^{\op} \to \Fun(I, \vect)$ is called \emph{flat} if the unit natural transformation $\eta_a \colon K(a, -) \to \mathcal{RL} K(a, -)$ is an isomorphism for every $a$ in $J$ for which $\parameterization(a) \neq 0$.
\end{step}

Every flat functor is thin. Since both left and right adjoints commute with direct sums, if $\parameterization$ is thin, then, for every free functor $C_0$ in $\Fun(J, \vect)$, the unit natural transformation $\eta_{C_0} \colon C_0 \to \mathcal{RL} C_0$ is also an epimorphism.
If $\parameterization$ is flat and $C_0$ in $\Fun(J, \vect)$ is a free functor, then the unit natural transformation $\eta_{C_0} \colon C_0 \to \mathcal{RL} C_0$ is an isomorphism if, and only if, $\parameterization(a) \neq 0$ for every $a$ in $\text{supp}(\beta^0 C_0)$.

\begin{prop}
\label{P:ThinRestriction}
Let $I$ and $J$ be finite posets.
If the functor $\parameterization \colon J^\op \to \Fun(I, \vect)$ is 
thin (resp.\ flat), then, for every subposet $L\subseteq J$, the restriction of $\parameterization$ to $L^\op \subseteq J^{\op}$ is also thin (resp.\ flat).
\end{prop}
\begin{proof}
For all $a$ in $J$, the functors $\mathcal{RL}K(a, -)$, $\mathcal{R}\parameterization(a)$, and $\Nat_I(\parameterization, \parameterization(a))$ are isomorphic. Thus $\parameterization$ is thin if, and only if, both of the following conditions are satisfied:
\begin{itemize}
\item for $a \not\leqJ b$ in $J$, $\Nat_I(\parameterization(b), \parameterization(a)) = 0$,
\item for $a \leqJ b$ in $J$, every natural transformation $\parameterization(b) \to \parameterization(a)$ is of the form
$\lambda \parameterization(a \leqJ b)$ for some $\lambda$ in the field $K$.
\end{itemize}
This characterization implies that if $\parameterization$ is thin, then its
restriction to $L^{\op}$ is also thin.

Similarly, the functor $\parameterization$ is flat if, and only if, both of the following conditions are satisfied:
\begin{itemize}
\item for $a \not\leqJ b$ in $J$, $\Nat_I(\parameterization(b), \parameterization(a)) = 0$,
\item for $a \leqJ b$ in $J$, if $\parameterization(a)\neq 0$, then $\Nat_I(\parameterization(b), \parameterization(a))$ is $1$-dimensional
and $\parameterization(a \leqJ b)$ is nonzero.
\end{itemize}
As before, this characterization implies that if $\parameterization$ is flat, then its restriction to $L^{\op}$ is also flat.
\end{proof}

Thinness is important because it implies that all the 
elements 
in the collection  $\mathcal{P}$ are
indecomposable (as their endomorphism algebras are 1 dimensional), and consequently this collection is independent (see~\ref{Rel:Freeness}).
Nevertheless we find the following proof  insightful.

\begin{prop}
\label{P:ThinIndependence}
Let $I$ and $J$ be finite posets.
If the functor $\parameterization \colon J^{\op} \to \Fun(I, \vect)$ is 
thin, then the collection $\mathcal{P} = \{\parameterization(a) \mid a\in J,\ \parameterization(a) \neq 0\}$ is independent (see~\ref{Rel:Freeness}). 
\end{prop}
\begin{proof}
Choose $a$ in $J$. The functor $\mathcal{RL} K(a, -) = \Nat_I(\parameterization, \parameterization(a)) = \mathcal{R}\parameterization(a)$ is the zero functor if, and only if, $\parameterization(a)$ is the zero functor. Thus, if $\parameterization(a)$ is nonzero, then the surjectivity of $\eta_a \colon K(a, -) \to \mathcal{RL} K(a, -)$ ($\parameterization$ is assumed to be thin) implies that $\eta_a$ is a minimal cover in $\Fun(J, \vect)$, in which case the standard $0^{\th}$ Betti diagram $\beta^0 \mathcal{R}\parameterization(a) = \beta^0 \mathcal{RL} K(a, -) \colon J \to \mathbb{N}$ has the following values:
\[
\beta^0 \mathcal{R}\parameterization(a)(b) = \begin{cases}
1 & \text{if } a = b, \\
0 & \text{otherwise.}
\end{cases}
\]

Let $\beta\colon J\to\mathbb{N}$ be a function whose support is finite and such that, if $\beta(a)\not=0$, then $\parameterization(a)\not=0$. Consider the functor $\bigoplus_{a\in J} \parameterization(a)^{\beta(a)}$, which is a finite sum because $J$ is finite.  Since $\mathcal R$ is right adjoint, it commutes with finite direct sums and consequently $\mathcal{R}\left(\bigoplus_{a\in J} \parameterization(a)^{\beta(a)}\right)$ is isomorphic to $\bigoplus_{a\in J} \mathcal{R}\parameterization(a)^{\beta(a)}$.
According to the calculation above, $\beta$ is the $0^\th$ Betti diagram of 
$\mathcal{R}\left(\bigoplus_{a\in J} \parameterization(a)^{\beta(a)}\right)$. The function  $\beta$ is therefore  determined uniquely by the isomorphism type of the functor $\bigoplus_{a\in J} \parameterization(a)^{\beta(a)}$, which gives the independence of $\parameterization$.
\end{proof}

If $\parameterization \colon J^{\op} \to \Fun(I,\vect)$ is thin, then we also have an effective way of constructing (minimal) $\mathcal{P}$-covers in $\Fun(I, \vect)$ using standard (minimal) covers in $\Fun(J,\vect)$. 

\begin{thm}
\label{T:MinCovers}
Let $I$ and $J$ be finite posets and $\parameterization \colon J^{\op} \to \Fun(I, \vect)$ a thin functor (see~\ref{D:Thinness}).
Then, for every $M \colon I \to \vect$, a natural transformation $C_0 \to \mathcal{R} M$ is a cover in $\Fun(J,\vect)$ if, and only if, its adjoint $\mathcal{L} C_0 \to M$ is a $\mathcal{P}$-cover in $\Fun(I, \vect)$.
Moreover, if $C_0 \to \mathcal{R} M$ is a minimal cover, then its adjoint $\mathcal{L} C_0 \to M$ is a minimal $\mathcal{P}$-cover. 
\end{thm}
\begin{proof}
Since, for every $a$ in $J$, the functor $\mathcal{L}K(a, -)$ is isomorphic to $\parameterization(a)$ and $\mathcal{L}$ commutes with direct sums, a functor $C_0$ in $\Fun(J, \vect)$ is free if, and only if, the functor $\mathcal{L} C_0$ in $\Fun(I, \vect)$ is $\mathcal{P}$-free. Recall that all projectives in $\Fun(J, \vect)$ are free.

Choose a functor $M \colon I \to \vect$, and a natural transformation $f \colon C_0 \to \mathcal{R}M$ with $C_0$ a free functor in $\Fun(J, \vect)$. Let $g \colon \mathcal{L} C_0 \to M$ be the adjoint of $f$. The following commuting diagram describes the relation between $f$ and $g$:
\[
\begin{tikzcd}
C_0
 \ar{r}{\eta_{C_0}}
 \ar[bend right, swap]{dr}{f}
&
\mathcal{RL} C_0
 \ar{d}{\mathcal{R} g}
\\
&
\mathcal{R} M
\end{tikzcd}
\]
By thinness assumption, the unit $\eta_{C_0} \colon C_0 \to \mathcal{RL} C_0$ is an epimorphism. Therefore $f$ is an epimorphism if, and only if, $\mathcal{R} g$ is an epimorphism, which, by Proposition~\ref{P:RelAndStdExactness}, happens if, and only if, $g$ is a $\mathcal{P}$-epimorphism. Thus $f$ is a cover if, and only if, $g$ is a $\mathcal{P}$-cover.

Next we discuss minimality. Suppose that $f \colon C_0 \to \mathcal{R} M$ is a minimal cover. Since we already showed that $g \colon \mathcal{L} C_0 \to M$ is a $\mathcal{P}$-cover, it remains to show its minimality. Let $h \colon \mathcal{L} C_0 \to \mathcal{L} C_0$ be an endomorphism of $g$. This endomorphism fits into the following commutative diagram:
\[
\begin{tikzcd}[row sep=5pt, column sep =50pt]
C_0
 \ar{r}[description]{\eta_{C_0}}
 \ar[dashed]{dd}
 \ar[bend left =30]{rrd}[description]{f}
&
\mathcal{RL} C_0
 \ar{dd}[swap]{\mathcal{R} h}
 \ar{dr}[description]{\mathcal{R} g}
\\
&&
\mathcal{R} M
\\
C_0
 \ar{r}[description]{\eta_{C_0}}
 \ar[bend right = 30]{rru}[description]{f}
&
\mathcal{RL} C_0
 \ar{ur}[description]{\mathcal{R} g}
\end{tikzcd}
\]
As before, the unit $\eta_{C_0} \colon C_0 \to \mathcal{RL} C_0$ is an epimorphism. The dashed arrow exists because $C_0$ is free. By the minimality of $f$, this dashed arrow is an isomorphism. The natural transformation $\mathcal{R} h$ is therefore an epimorphism.
By Corollary~\ref{C:RelAndStdEpi}, $h \colon \mathcal{L} C_0 \to \mathcal{L} C_0$ is a $\mathcal{P}$-epimorphism. Since $\mathcal{L} C_0$ is $\mathcal{P}$-free, $h$ being a $\mathcal{P}$-epimorphism means that $h$ is an epimorphism.
As the values of $\mathcal{L}C_0$ are finite dimensional vector spaces, $h$ is an isomorphism.
The natural transformation $g$ is therefore a minimal cover.
\end{proof}

\begin{example}
Considering the functor from Example~\ref{E:M}
\[
M =
\begin{tikzpicture}[scale = \diagramscale, vcenter=1]
\tikzgrid{1}
\fill[functor fill, functor draw] (0, 0) rectangle (1, 1);
\end{tikzpicture}\;,
\]
the minimal cover of $\mathcal{R} M$ is
\[
\begin{tikzpicture}[scale = \diagramscale, vcenter=2]
\tikzgrid{2}
\fill[functor fill, functor draw] (0, 1) rectangle (3, 3);
\fill[functor fill, functor draw] (1, 0) rectangle (3, 3);
\end{tikzpicture}
\; \longrightarrow \;
\begin{tikzpicture}[scale = \diagramscale, vcenter=2]
\tikzgrid{2}
\fill[functor fill, functor draw] (0, 1) -- (1, 1) -- (1, 0) -- (2, 0) -- (2, 2) -- (0, 2) -- cycle;
\end{tikzpicture}\;,
\]
that is, $K((00, 10), -) \oplus K((00, 01), -) \rightarrow \mathcal{R} M$. Then, by Theorem \ref{T:MinCovers}, the minimal $\mathcal{P}$-cover of $M$ is obtained by applying $\mathcal{L}$ (see the discussion before \ref{E:M}):
\[
\begin{tikzpicture}[scale = \diagramscale, vcenter=1]
\tikzgrid{1}
\fill[functor fill, functor draw] (0, 0) rectangle (2, 1);
\fill[functor fill, functor draw] (0, 0) rectangle (1, 2);
\end{tikzpicture}
\; \longrightarrow \;
\begin{tikzpicture}[scale = \diagramscale, vcenter=1]
\tikzgrid{1}
\fill[functor fill, functor draw] (0, 0) rectangle (1, 1);
\end{tikzpicture}\;.
\]
\end{example}

Theorem~\ref{T:MinCovers} has several important consequences, including the acyclicity of the collection $\mathcal{P}$:

\begin{cor}
\label{P:ThinAcyclicity}
Let $I$ and $J$ be finite posets. Suppose the functor $\parameterization \colon J^{\op} \to \Fun(I, \vect)$ is thin. Then
\begin{enumerate}
\item every functor in $\Fun(I, \vect)$ admits a minimal $\mathcal{P}$-resolution (see~\ref{Rel:Resolutions});
\item the collection $\mathcal{P} = \{\parameterization(a) \mid a \in J,\ \parameterization(a) \neq 0\}$ is acyclic (see~\ref{Rel:Projectives}). 
\end{enumerate}
\end{cor}
\begin{proof}
\noindent
(1) Every functor in $\Fun(J, \vect)$ admits a minimal cover (see~\ref{Std:MinCovers}). 
Thus according to Theorem~\ref{T:MinCovers} every functor in $\Fun(I, \vect)$ admits a minimal $\mathcal{P}$-cover, and hence a minimal $\mathcal{P}$-resolution.
\smallskip

\noindent
(2) Let $M$ be a $\mathcal{P}$-projective object in $\Fun(I, \vect)$. We need to show that it is $\mathcal{P}$-free.
Consider a minimal cover $C_0 \to \mathcal{R} M$ in $\Fun(J, \vect)$ (see~\ref{Std:MinCovers}).
Its adjoint $\mathcal{L} C_0 \to M$, according to Theorem~\ref{T:MinCovers}, is then a minimal $\mathcal{P}$-cover of $M$ in $\Fun(I, \vect)$, and therefore has to be an isomorphism since $M$ is $\mathcal{P}$-projective.
Note that $\mathcal{L} C_0$ is $\mathcal{P}$-free, and consequently so is $M$.
\end{proof}

To summarize, if the functor $\parameterization \colon J^{\op} \to \Fun(I, \vect)$ is thin, then the collection $\mathcal{P}$ is independent (see~\ref{P:ThinIndependence}) and acyclic (see~\ref{P:ThinAcyclicity}.(2)), and every functor $M \colon I \to \vect$ has a minimal $\mathcal{P}$-resolution (see~\ref{P:ThinAcyclicity}.(1)).
Thus thinness guarantees that the $\mathcal{P}$-Betti diagram $\beta^d_{\mathcal{P}} M \colon \mathcal{P} \to \mathbb{N}$ is well defined for every $d \geq 0$ (see~\ref{Rel:Betti}). The rest of this section is devoted to describing methods to calculate some values of these $\mathcal{P}$-Betti diagrams. We start by looking at what we can prove for a general poset $J$:

\begin{prop}
\label{P:RelativeZerothBetti}
Let $I$, $J$ be finite posets, and $\parameterization \colon J^\op \to \Fun(I, \vect)$ and $M \colon I\to \vect$ be functors.
\begin{enumerate}
\item If $\parameterization$ is thin, then $\beta^0_{\mathcal{P}} M(\parameterization(a)) = \beta^0(\mathcal{R} M)(a)$ for every $a$ in $J$ for which $\parameterization(a) \neq 0$.
\item If $\parameterization$ is flat, then $\beta^d_{\mathcal{P}} M(\parameterization(a)) = \beta^d(\mathcal{R} M)(a)$ for every $d\geq 0$ and every $a$ in $J$ for which $\parameterization(a) \neq 0$.
\end{enumerate}
\end{prop}
\begin{proof}
Let $f\colon C_0 \to \mathcal{R} M$ be a minimal cover. By definition of $\beta^0(\mathcal{R} M)$, the functor $C_0$ is isomorphic to $\bigoplus_{a \in J} K(a, -)^{\beta^0(\mathcal{R} M)(a)}$.
Recall that, for every $a$ in $J$, the functor $\mathcal{L}K(a,-)$ is isomorphic to $\parameterization(a)$.
This, together with the fact that $\mathcal{L}$ commutes with direct sums, implies that $\mathcal{L} C_0$ is isomorphic to the functor
\[
\bigoplus_{\substack{a \in J \\ \parameterization(a) \neq 0}} \parameterization(a)^{\beta^0(\mathcal{R} M)(a)}.
\]

Suppose $\parameterization$ is thin. According to Theorem~\ref{T:MinCovers}, the minimal $\mathcal{P}$-cover of $M$ is given by the adjoint $g\colon \mathcal{L} C_0\to M$ of the minimal cover $f\colon C_0 \to \mathcal{R} M$. We conclude that the equality $\beta^0_{\mathcal{P}} M(\parameterization(a)) = \beta^0(\mathcal{R} M)(a)$ holds for all $a$ for which $\parameterization(a) \neq 0$. This shows (1).

Suppose $\parameterization$ is flat. 
The natural transformations
$\ker(g) \hookrightarrow \mathcal{L} C_0$ and $g \colon \mathcal{L} C_0 \to M$ form a $\mathcal{P}$-exact sequence in $\Fun(I, \vect)$. According to Proposition~\ref{P:RelAndStdExactness}, after applying $\mathcal{R}$, the natural transformations $\mathcal{R}(\ker(g) \hookrightarrow \mathcal{L} C_0)$ and $\mathcal{R} g$ form therefore an exact sequence in $\Fun(J, \vect)$. Since $\mathcal{R}$ is right adjoint, it preserves monomorphisms and hence $\mathcal{R} \ker(g)$ is the kernel of $\mathcal{R} g$. 
We can then form the following commutative diagram with the indicated arrows representing epimorphisms and monomorphisms, and with the vertical sequences being exact:
\[
\begin{tikzcd}
\ker(f)
  \ar[hook]{d}
  \ar{r}
&
\mathcal{R} \ker(g)
  \ar[hook]{d} 
\\
C_0
  \ar[two heads]{d}[swap]{f}
  \ar{r}{\eta}
&
\mathcal{RL} C_0
  \ar[two heads]{d}{\mathcal{R} g}
\\
\mathcal{R}M
  \ar[equal]{r}
&
\mathcal{R} M
\end{tikzcd}
\]
By flatness, $\mathcal{R}\mathcal{L}C_0$ is free and hence, by minimality of $f$, the unit $\eta$ is an isomorphism. The natural transformation represented by the top horizontal arrow in this diagram is then also an isomorphism.
Consequently, for every $d \geq 0$ and every $a$ in $J$, $\beta^d(\ker(f))(a) = \beta^d(\mathcal{R} \ker(g))(a)$, and we conclude with the following sequence of equalities:
\begin{align*}
\beta^d_{\mathcal{P}} M(\parameterization(a)) & = \beta^{d - 1}_\mathcal{P} (\ker(g))(\parameterization(a)) \tag{Since $g \colon \mathcal{L} C_0 \to M$ is a minimal $\mathcal{P}$-cover} \\
& = \beta^{d - 1} (\mathcal{R} \ker(g))(a) \tag{By induction} \\
& = \beta^{d - 1} (\ker(f))(a) \tag{By isomorphism} \\
& = \beta^d (\mathcal{R} M)(a). \tag{Since $f \colon C_0 \to \mathcal{R} M$ is a minimal cover}
\end{align*}
This proves (2).
\end{proof}

Proposition~\ref{P:RelativeZerothBetti}.(1) provides an algorithm to calculate the $0^\th$ $\mathcal{P}$-Betti diagram of a functor $M \colon I \to \vect$ when $\parameterization$ is thin: first take $\mathcal{R} M \colon J \to \vect$, then calculate its standard $0^\th$ Betti diagram (over the poset $J$), and finally restrict the obtained $0^\th$ Betti diagram to the $a$ in $J$ for which $\parameterization(a) \neq 0$.
The standard $0^\th$ Betti diagrams can be calculated for example using radicals (see~\ref{Std:MinCovers}).
This algorithm can be then used inductively in order to calculate the $\mathcal{P}$-Betti diagrams $\beta^d_{\mathcal{P}} M$ of $M$ for all $d \geq 0$, as explained in~\ref{Rel:Betti}. 
In every step of this procedure we need to evaluate 
the functor $\mathcal{R}$ on successive kernels (see the second step in the sequence of equalities at the end of the proof of~\ref{P:RelativeZerothBetti}). 
We would like to avoid that step, for example by showing beforehand that
the $\mathcal{P}$-Betti diagrams over $I$ at $\parameterization(a)$ are equal to the corresponding Betti diagrams over $J$ at $a$.
This is the case if, for example,
$\parameterization$ is flat (see Proposition~\ref{P:RelativeZerothBetti}.(2)). However, our key examples (see Section~\ref{S:Examples}) are not flat but thin.
For an arbitrary thin functor 
there could be elements $a$ in $J$ for which the numbers $\beta^d_{\mathcal{P}} M(\parameterization(a)) \neq \beta^d (\mathcal{R} M)(a)$ differ. Such elements are called $\parameterization$-degenerate:

\begin{step}[Degeneracy locus]
\label{D:DegeneracyLocus}
Let $\parameterization \colon J^{\op} \to \Fun(I, \vect)$ be thin.
An element $a$ in $J$ is called \emph{$\parameterization$-degenerate} if $\parameterization(a) = 0$ or if $\parameterization(a) \neq 0$ and there exists a functor $M \colon I \to \vect$ and a natural number $d \geq 0$, for which $\beta^d_{\mathcal{P}} M(\parameterization(a)) \neq \beta^d (\mathcal{R}M)(a)$. 
The collection of $\parameterization$-degenerate elements is called \textbf{degeneracy locus} of $\parameterization$. 
\end{step}

If $\parameterization$ is flat, then its degeneracy locus is given by $\{a \in J \mid \parameterization(a) = 0\}$ (see~\ref{P:RelativeZerothBetti}.(2)). For an arbitrary thin functor $\parameterization$ we do not have an explicit description of its degeneracy locus. However, we can approximate it when $J$ is an upper semilattice.

\begin{thm}[Main result, first half]
\label{T:ContainDegeneracy}
Let $I$ be a finite poset, $J$ a finite upper semilattice, and $\parameterization \colon J^\op \to \Fun(I, \vect)$ a thin functor. Then the degeneracy locus of $\parameterization$ (see~\ref{D:DegeneracyLocus}) is contained in the set
\[
\bigcup_{\substack{a \in J \\ d \geq 0}} \supp(\beta^d \ker \eta_a) \subseteq J.
\]
\end{thm}
\begin{proof}
Write $\mathcal{G} := \bigcup_{a, d} \supp(\beta^d \ker \eta_a)$. Let us prove the contraposition, that for $M \colon I \to \vect$ a functor, $a$ in $J \setminus \mathcal{G}$, and $d \geq 0$, we have $\parameterization(a) \neq 0$ and $\beta^d_\mathcal{P} M(\parameterization(a)) = \beta^d (\mathcal{R} M)(a)$. Note that if $\parameterization(a) = 0$, then $\supp(\beta^0 \ker \eta_a) = \{a\}$, so the first statement is true.
We proceed by induction on $d \geq 0$ for the second statement. The case $d = 0$ is the statement of Proposition~\ref{P:RelativeZerothBetti}.
Now let $d > 0$.

We proceed exactly as in the proof of Proposition~\ref{P:RelativeZerothBetti}.
Let $M \colon I \to \vect$ be a functor and consider the following commutative diagram:
\[\begin{tikzcd}
\ker(f)
  \ar[hook]{d}
  \ar{r}
&
\mathcal{R} \ker(g)
  \ar[hook]{d} 
\\
C_0
  \ar[two heads]{d}[swap]{f}
  \ar{r}{\eta}
&
\mathcal{RL} C_0
  \ar[two heads]{d}{\mathcal{R} g}
\\
\mathcal{R} M
  \ar[equal]{r}
&
\mathcal{R} M
\end{tikzcd}\]
where $f$ is a minimal cover in $\Fun(J, \vect)$, $g$ is its left adjoint, which
according to Theorem~\ref{T:MinCovers} is a minimal $\mathcal{P}$-cover in $\Fun(I, \vect)$, the indicated arrows are epimorphisms and monomorphisms, and the vertical sequences are exact.
The thinness assumption implies that the unit natural transformation $\eta$ (the middle horizontal arrow in the above diagram) is an epimorphism.
This, together with the exactness of the vertical sequences, implies that the natural transformation represented by the top horizontal arrow is also an epimorphism.
Also by exactness of the vertical sequences, the kernel of this natural transformation is isomorphic to $\ker(\eta)$. In this way we obtain an exact sequence in $\Fun(J, \vect)$ of the form
\[
\begin{tikzcd}
0
  \ar{r}
&
\ker(\eta)
  \ar{r}
&
\ker(f)
  \ar{r}
&
\mathcal{R} \ker(g)
  \ar{r}
& 0.
\end{tikzcd}
\]
Let us show that for all $a$ in $J \setminus \mathcal{G}$ and $d \geq 0$, $\beta^d(\ker(f))(a) = \beta^d(\mathcal{R} \ker(g))(a)$. Since $C_0$ is free and both $\mathcal{R}$ and $\mathcal{L}$ commute with direct sums, $\ker(\eta)$ is isomorphic to a direct sum of functors of the form $\ker(\eta_a\colon K(a,-) \to \mathcal{RL} K(a, -))$. By definition of $\mathcal{G}$, we have the containment $\supp(\beta^d \ker(\eta)) \subseteq \mathcal{G}$ for every $d \geq 0$. Thus, for every $d \geq 0$ and $a$ in $J \setminus \mathcal{G}$, $\beta^d(\ker(\eta))(a) = 0$. Since $J$ is an upper semilattice by assumption, we apply Corollary~\ref{C:EqualBetti} to deduce that $\beta^d(\ker(f))(a) = \beta^d(\mathcal{R} \ker(g))(a)$.

We conclude with the following sequence of equalities:
\begin{align*}
\beta^d_{\mathcal{P}} M(\parameterization(a)) & = \beta^{d - 1}_\mathcal{P} (\ker(g))(\parameterization(a)) \tag{Since $g \colon \mathcal{L} C_0 \to M$ is a minimal $\mathcal{P}$-cover} \\
& = \beta^{d - 1} (\mathcal{R} \ker(g))(a) \tag{By induction} \\
& = \beta^{d - 1} (\ker(f))(a) \tag{By our previous claim} \\
& = \beta^d (\mathcal{R} M)(a). \tag{Since $f \colon C_0 \to \mathcal{R} M$ is a minimal cover}
\end{align*}
\end{proof}

With some additional hypotheses on the parameterization, we get an exact description of the degeneracy locus.

\begin{cor}[Main result, second half]
\label{C:BettiFromDegeneracy}
Let $I$ be a finite poset, $J$ a finite upper semilattice, and $\parameterization \colon J^\op \to \Fun(I, \vect)$ a thin functor. Suppose that, for all $a$ in $J$, the sublattice $\langle \supp(\beta^0 \ker \eta_a)\rangle$ is contained in $\{b \in J \mid \parameterization(b) = 0\}$.
Then the degeneracy locus of $\parameterization$ is 
$\{b\in J \mid \parameterization(b)=0\}$, which is equal to $\cup_{a \in J} \langle \supp(\beta^0 \ker \eta_a) \rangle$.

Explicitly, for every functor $M \colon I \to \vect$, every element $a$ in $J$ such that $\parameterization(a) \neq 0$, and every $d \geq 0$,
\[
\beta^d_{\mathcal{P}} M(\parameterization(a)) = \dim H_d(\mathcal{K}_a \mathcal{R} M).
\]
\end{cor}

\begin{proof}
For $a$ in $J$ and $d \geq 0$, $\supp(\beta^d \ker \eta_a)$ is contained in $\langle \supp(\beta^0 \ker \eta_a) \rangle$, which itself is contained in $\{b \in J \mid \parameterization(b) = 0\}$ by hypothesis. By Theorem~\ref{T:ContainDegeneracy}, we deduce that the degeneracy locus of $\parameterization$ is also contained in $\{b \in J \mid \parameterization(b) = 0\}$. In fact, they are equal, and they are both equal to $\cup_{a \in J} \langle \supp(\beta^0 \ker \eta_a) \rangle$. By definition of the degeneracy locus and Theorem~\ref{T:BettiKoszul} respectively, for all $M \colon I \to \vect$, $a$ in $J$ such that $\parameterization(a) \neq 0$, and $d \geq 0$, we get the equalities
\[
\beta^d_{\mathcal{P}} M(\parameterization(a)) = \beta^d \mathcal{R} M(a) = \dim H_d(\mathcal{K}_a \mathcal{R} M).
\qedhere
\]
\end{proof}

\begin{example}
Our running Example~\ref{E:Running} satisfies the conditions of Corollary \ref{C:BettiFromDegeneracy}. This is proved in general for lower hooks in \ref{E:LowerHooks}, but here we can check this by hand: for instance, for $a = (00, 10)$, we computed the unit natural transformation $\eta_a$ in \ref{E:ManualThinness}. We find that
\[
\ker \eta_a = 
\begin{tikzpicture}[scale = \diagramscale, vcenter=2]
\tikzgrid{2}
\fill[functor fill, functor draw] (2, 0) rectangle (3, 3);
\end{tikzpicture}\;,
\]
and so $\supp(\beta^0 \ker \eta_a) = \{(10, 10)\} = \langle \supp(\beta^0 \ker \eta_a) \rangle$, and from \ref{E:Running} we know that $\parameterization((10, 10)) = 0$. We reason similarly for all other $a$ in $J$.
\end{example}

\section{Examples}
\label{S:Examples}
This last section is devoted to
examples of independent and acyclic collections
together with  calculations of the induced relative Betti diagrams in the category
$\Fun(I, \vect)$, where $(I, \leq)$ is a finite poset. Our strategy for performing these calculations is to find an appropriate parameterization of the collection  via an upper semilattice $J$ and 
 a functor $\parameterization\colon J^\op\to \Fun(I, \vect)$ as explained in Section~\ref{S:GradedRelHomAlg}.
 
\halfsection{Singleton}
\label{S:Singleton}
We start with the simplest nonempty example.
Choose a functor $P_0\colon I\to \vect$ and consider the collection $\{P_0\}$ consisting only of this  chosen functor. The simplest way to parameterize this collection is to consider the poset $[0]$  with a single element $0$ and the functor $\parameterization \colon [0]^\op \to \Fun(I, \vect)$ sending $0$ to $P_0$.

We do a case analysis by the dimension of $\Nat_I(P_0, P_0)$.
If $\dim \Nat_I(P_0, P_0) = 0$, then $P_0 = 0$ and the collection $\mathcal{P}=\{\parameterization(a)\mid a\in [0],\ \parameterization(a)\not=0\}$ is empty. For every functor $M \colon I \to \vect$, $0 \to M$ is a minimal $\mathcal{P}$-resolution, and the Betti diagrams of $M$ are therefore given by the empty function $\varnothing \to \mathbb{N}$.

If $\dim \Nat_I(P_0, P_0) = 1$, then the functor $\parameterization$ is flat (explicitly, $K(0, -) \to \Nat_I(P_0, P_0)$ is an isomorphism).  For every functor $M \colon I \to \vect$, the minimal cover of $\mathcal{R} M = \Nat_I(P_0, M)$ is the identity, so by Theorem~\ref{T:MinCovers}, the minimal $\mathcal{P}$-cover of $M$ is its left adjoint, which is the evaluation natural transformation $P_0 \otimes_K \Nat_I(P_0, M) \to M$.
By Proposition~\ref{P:RelativeZerothBetti}:
\[
\beta^d_{\mathcal{P}}M(P_0)=
\begin{cases}
\dim \Nat_I(P_0, M) & \text{ if } d = 0, \\
0 & \text{ if } d > 0.
\end{cases}
\]

If $\dim \Nat_I(P_0, P_0) > 1$, then the functor $\parameterization$ is not thin, and we observe for instance that the minimal cover of $\mathcal{R} P_0$ is adjoint to the $\mathcal{P}$-cover $P_0 \otimes_K \Nat_I(P_0, P_0) \to P_0$, which is not the minimal $\mathcal{P}$-cover of $P_0$.

\bigskip

We now proceed with the remaining examples, studying homological algebra relative to various collections of functors in $\Fun(I, \vect)$. In each example, to the extent that is possible, we proceed as follows.
We start by defining a \textbf{collection} and its \textbf{parameterization} by a functor $\parameterization \colon J^\op \to \Fun(I, \vect)$. We then check if the functor $\parameterization$ is flat or thin.
\begin{itemize}
\item If $\parameterization$ is \textbf{flat}, then we can apply Proposition~\ref{P:RelativeZerothBetti}.(2).
\item If $\parameterization$ is \textbf{thin}, then our aim is to apply Corollary~\ref{C:BettiFromDegeneracy}. To do this we need to verify the other conditions of the corollary, namely that the poset $J$ is an \textbf{upper semilattice} and that $\parameterization$ sends the sublattice $\langle \supp(\beta^0 \ker \eta_a) \rangle$ to $0$ for all $a$ in $J$. We call this latter assumption the \textbf{degeneracy} condition.
\end{itemize}
In both cases, we conclude that we can compute Betti diagrams relative to $\mathcal{P}$ from Koszul complexes over $J$. To do this effectively, we characterize \textbf{parents} in $J$ and then give an explicit formula for \textbf{Koszul complexes}. Then, for every functor $M \colon I \to \vect$, element $a$ in $J$ such that $\parameterization(a) \neq 0$, and $d \geq 0$, we have the equality
\[
\beta^d_{\mathcal{P}} M(\parameterization(a)) = \dim H_d(\mathcal{K}_{a} \mathcal{R} M).
\]
Finally, we compute an explicit relative projective resolution of the following functor $M_0$ on a $6 \times 6$ grid, whose minimal standard resolution is illustrated below:
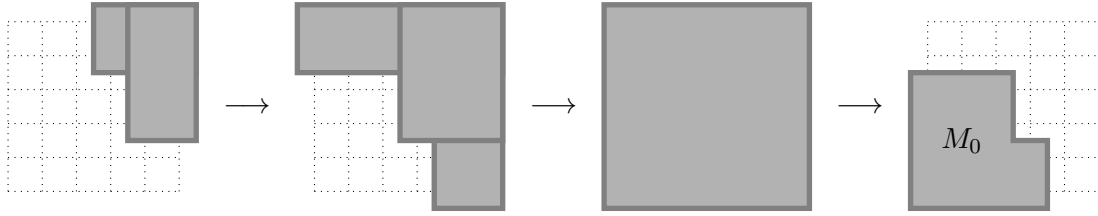
\[
\begin{tikzpicture}[scale = \diagramscale, vcenter=5]
\tikzgrid{5}
\fill[functor fill, functor draw] (3, 4) rectangle (6, 6);
\fill[functor fill, functor draw] (4, 2) rectangle (6, 6);
\end{tikzpicture}
\;\; \longrightarrow \;\;
\begin{tikzpicture}[scale = \diagramscale, vcenter=5]
\tikzgrid{5}
\filldraw[functor fill, functor draw] (0, 4) rectangle (6, 6);
\filldraw[functor fill, functor draw] (4, 0) rectangle (6, 6);
\filldraw[functor fill, functor draw] (3, 2) rectangle (6, 6);
\end{tikzpicture}
\;\; \longrightarrow \;\;
\begin{tikzpicture}[scale = \diagramscale, vcenter=5]
\tikzgrid{5}
\filldraw[functor fill, functor draw] (0, 0) rectangle (6, 6);
\end{tikzpicture}
\;\; \longrightarrow \;\;
\begin{tikzpicture}[scale = \diagramscale, vcenter=5]
\tikzgrid{5}
\filldraw[functor fill, functor draw] (0, 0) -- (0, 4) -- (3, 4) -- (3, 2) -- (4, 2) -- (4, 0) -- cycle;
\node at (1.5, 2) {$M_0$};
\end{tikzpicture}
\]
As in Example~\ref{E:Running}, the functor $M_0$ illustrated here is equal to $K$ on its support (vertices in the shaded area) and $0$ elsewhere, with all the transition
functions $M_0(v\leq w)$ having maximal rank.

\halfsection{All subfunctors}
\label{S:E:Subfunctors}

\subsection*{Collection}
We consider the collection of all nonzero elements in 
$\text{Sub}(I)$, the set of all subfunctors of the constant functor $K_I$ (see Section~\ref{S:FiltrationsSubfunctors}).

\subsection*{Parameterization}
Recall
that there is a poset relation $\leqJ$ on $\text{Sub}(I)$
(see~\ref{D:Antichains})
such that $F \leqJ G$ 
if, and only if, $\supp(F) \supseteq \supp(G)$, which is  equivalent to 
$\supp(F) \supseteq \supp(\beta^0 G)$. 
According to~\ref{D:SubfunctorsUpsets} and~\ref{D:Antichains},
the poset $(\text{Sub}(I), \leqJ)$ can be identified with either 
 the opposite of the inclusion poset
$(\text{Up}(I),\subseteq)$ of upsets of $I$, or with 
the poset   $(\text{Anti}(I),\leqJ)$ of antichains in $I$. In particular $(\text{Sub}(I), \leqJ)$ is an upper semilattice.

Define $\parameterization \colon (\text{Sub}(I), \leqJ)^\op \to \Fun(I, \vect)$ to be the identity on objects, and to map a relation $F \leqJ G$ to the corresponding natural transformation $G \subseteq F$.
The collection $\mathcal{P}=\{\parameterization(F)\mid F\in\text{Sub}(I),\ \parameterization(F)=F\not=0\}$ consists of all nonzero subfunctors of $K_I$.

\subsection*{Flatness}
This requires the additional assumption of $I$ having a unique maximal element.
\begin{prop}
Let $I$ be a finite poset with a unique maximal element. Then the parameterization $\parameterization$ is flat.
\end{prop}
\begin{proof}
Let $F, G \subseteq K_I$ be nonzero subfunctors. By \cite[Prop.~3.10(2)]{Miller2020}, the dimension of the space $\Nat(G, F)$ is equal to the number of connected components of $\supp(G)$ fully contained in $\supp(F)$. By hypothesis of $I$ having a unique maximal element, $\supp(F)$ and $\supp(G)$ each have only one connected component, so $\dim(\Nat(G, F)) = 1$ if $F \leqJ G$ and $0$ otherwise. We deduce that the unit natural transformation $\eta_F \colon K(F, -) \to \mathcal{RL} K(F, -)$ is an isomorphism, and so we conclude that $\parameterization$ is flat.
\end{proof}

According to Proposition~\ref{P:RelativeZerothBetti}.(2), the flatness of $\parameterization$
gives the equality $\beta^d_{\mathcal{P}} M(F) = \beta^d(\mathcal{R} M)(F)$ for every functor $M \colon I \to \vect$, nonzero subfunctor $F \subseteq K_I$, and $d \geq 0$.
The $\mathcal{P}$-Betti diagrams over $I$ can be therefore expressed as the standard Betti diagrams over 
the upper semilattice
$(\text{Sub}(I),\leqJ)$. Thus these standard Betti diagrams can be calculated using Koszul complexes, once we have identified and enumerated parents and their meets in $(\text{Sub}(I), \leqJ)$
(see~\ref{C:BettiFromDegeneracy}).

\subsection*{Parents}
We identify $(\text{Sub}(I),\leqJ)$ with the opposite poset of upsets $\Up(I)^\op$ via the map $F \mapsto \supp(F)$ (see~\ref{D:SubfunctorsUpsets}).
For a subset $S$ of $I$, denote by $S^c$ its complement in $I$ and by $\Max(S)$ the set of maximal elements of $S$.
The sets of parents in $\Up(I)^\op$ can be described as follows:
\begin{lemma}
\label{L:UpsetParents}
The set of parents of  $U$ in  $\Up(I)^\op$
is  $\{U \cup \{v\} \mid v \in \Max (U^c)\}$.
\end{lemma}
\begin{proof}
Let $v$ be an element of $\Max (U^c)$. The set $U \cup \{v\}$ is indeed an upset because every element greater than $v$ 
is not in $U^c$, and hence it is in $U$. Since $U \cup \{v\}$ contains $U$, the former is indeed a lower bound of the latter in $\textnormal{Up}(I)^\op$. Finally, since there is only a one-element difference, $U \cup \{v\}$ is a parent of $U$.

Conversely, let $V$ be a parent of $U$ in $\textnormal{Up}(I)^\op$. Then $U$ is properly contained in $V$. Let $v$ be a maximal element in $V \setminus U$. Since $V$ is an upset, $v$ is also an element of $\Max (U^c)$, and consequent  $V=U \cup \{v\}$. 
\end{proof}

Recall that for an upset $U$ of $I$, the symbol $K_U$ denotes the unique subfunctor of $K_I$ for which $\supp(K_U)=U$. Thus, in terms of subfunctors, according to Lemma~\ref{L:UpsetParents}, the parents of $F\subseteq K_I$ 
in the poset $(\text{Sub}(I),\leqJ)$
are subfunctors of the form $K_{\supp(F) \cup \{v\}}\subseteq K_I$ for $v$ in $\Max(\supp(F)^c)$.

\subsection*{Koszul complexes}
For  $M \colon I \to \vect$ and
$F\subseteq K_I$, the Koszul complex of $\mathcal{R} M\colon \text{Sub}(I)\to \vect$ at $F$ is therefore given  in degree $d \geq 0$ by
\[
(\mathcal{K}_F \mathcal{R} M)_d = \bigoplus_{\substack{S \subseteq \Max (\text{supp}(F)^c) \\ |S| = d}} \Nat(K_{\text{supp}(F) \cup S}, M),
\]
with the differentials as defined  in \ref{Std:KoszulComplexes}.
The value $\beta^d_{\mathcal{P}} M(F)$ is then given by the dimension of  the $d^{\th}$ homology of this chain complex. For example, in the case $I$ is the  $6 \times 6$ grid,  the relative projective resolution of $M_0$ can be illustrated as
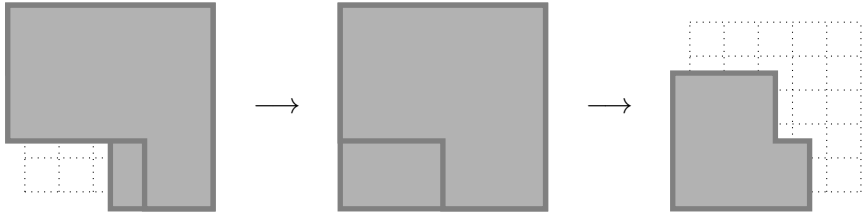
\[
\begin{tikzpicture}[scale = \diagramscale, vcenter=5]
\tikzgrid{5}
\filldraw[functor fill, functor draw] (6, 6) -- (0, 6) -- (0, 4) -- (3, 4) -- (3, 0) -- (6, 0) -- cycle;
\filldraw[functor fill, functor draw] (6, 6) -- (0, 6) -- (0, 2) -- (4, 2) -- (4, 0) -- (6, 0) -- cycle;
\end{tikzpicture}
\quad \longrightarrow \quad
\begin{tikzpicture}[scale = \diagramscale, vcenter=5]
\tikzgrid{5}
\filldraw[functor fill, functor draw] (0, 0) rectangle (6, 6);
\filldraw[functor fill, functor draw] (6, 6) -- (0, 6) -- (0, 2) -- (3, 2) -- (3, 0) -- (6, 0) -- cycle;
\end{tikzpicture}
\quad \longrightarrow \quad
\begin{tikzpicture}[scale = \diagramscale, vcenter=5]
\tikzgrid{5}
\filldraw[functor fill, functor draw] (0, 0) -- (0, 4) -- (3, 4) -- (3, 2) -- (4, 2) -- (4, 0) -- cycle;
\end{tikzpicture}
\]

Since $\parameterization \colon \text{Sub}(I)^\op \to \Fun(I, \vect)$ is flat, so are all of its restrictions (see~\ref{P:ThinRestriction}). 
Thus one way of producing examples of flat functors for which we can use Koszul complexes to calculate the associated relative Betti diagrams is to find subposets of $\text{Sub}(I)$ which are upper semilattices (not necessarily sublattices of $(\text{Sub}(I),\leqJ)$) and for which we can identify parents of its elements and meets of subsets of parents that are bounded below.
The next example is an illustration of this strategy.

\halfsection{Translated functors}
\label{S:E:Translated functors}

\subsection*{Collection}
Suppose the poset $(I, \leq)$ is the product $\{0 < \cdots < n\}^r$, which is a distributive lattice. Its elements will be denoted as words $w_1\cdots w_r$.
Fix an antichain $S$ in $\text{Anti}(I)$ and define $I - S := \{v \in I \mid S + v\subseteq I\}$.
The set $I - S$ consists of the elements $v$ in $I$ for which 
all the coordinates of $t + v$ are bounded by $n$ for every $s$ in $S$. Thus 
\[
I - S = \prod_{i = 1}^r \left\{0, \ldots, n - \max\{s_i \mid s \in S\}\right\}\ \subseteq\ \{0, \ldots, n\}^r = I. 
\]
Moreover, the subposet $(I - S) \subseteq I$ is a sublattice. 

For $v$ in $I - S$, the subset $S+v\subseteq I$ is also an antichain in $I$. The associated
subfunctor $K(S + v, -)\subseteq K_I$ is called the \emph{translation} of $K(S, -) \subseteq K_I$ by $v$
(see~\ref{D:Antichains} for the notation). 
In this example we consider the collection of 
all such translations 
$\mathcal{P} := \{K(S+v,-)\mid v \in I-S\}$.

\subsection*{Parameterization}
The lattice $I - S$ provides a natural choice for parameterizing this  collection because, for $v$ and $w$ in $I - S$, the inclusion $K(S + v, -) \subseteq K(S + w, -)$ exists if, and only if, $v \geq w$ in $I$. In particular, $K(S + v, -) = K(S + w, -)$ if, and only if, $v = w$.
Set $\parameterization \colon (I - S)^\op \to \Fun(I, \vect)$ to be the functor which maps $v$ in $(I - S)^\op$ to the translation $K(S + v, -)$.

\subsection*{Flatness}
This functor is flat by Proposition~\ref{P:ThinRestriction}, since it can be identified with the restriction of the functor $\parameterization$ discussed in~\ref{S:E:Subfunctors} to the subposet 
$\{K(S + v, -) \mid v \in I - S\} \subseteq \text{Sub}(I)$. Although the inclusion of posets may fail to be a sublattice inclusion, the poset $I - S$ is a lattice (since it is a sublattice of $I$) and consequently the $\mathcal{P}$-Betti diagrams  can be calculated via Koszul complexes.

\subsection*{Parents} 
For every element $v$ of $(I- S) \subseteq I$, the parents of $v$ in $I$ are also in $I - S$. Thus the parents of $v$ in $I - S$ are the same as those of $v$ in $I$.

\subsection*{Koszul complexes}
For $M \colon I \to \vect$ and  $v$ in $I - S$, the Koszul complex of 
$\mathcal{R} M\colon (I - S)\to \vect$ at $v$ is given in degree $d \geq 0$ by
\[
(\mathcal{K}_v \mathcal{R} M)_d = \bigoplus_{\substack{T \subseteq \mathcal{U}_{I-S}(v)\\ |T| = d }} \Nat(K({\textstyle S + \bigwedge_{(I-S)} T}, -), M),
\]
with differentials as defined in \ref{Std:KoszulComplexes}. 
 Finally, the value $\beta^d_{\mathcal{P}} M(K(S+v,-))$ is given by the dimension of  the $d^{\th}$ homology of this chain complex.
For example, if we choose 
$I=\{0<\cdots <5\}^2$ and 
$S=\{02, 10\}\subset I$, then the projective resolution of $M_0$ relative to translations of the subfunctor $K(\{02, 10\}, -)\subseteq K_I$ can be illustrated as
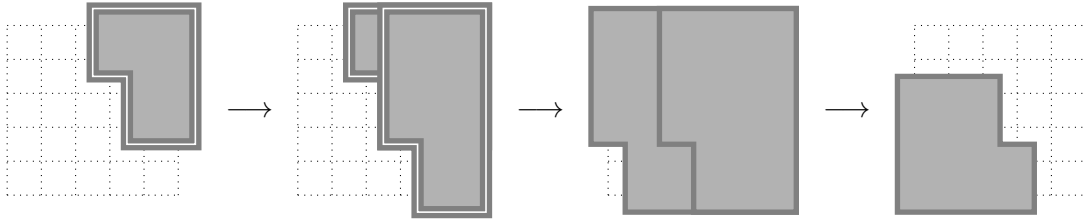
\[
\begin{tikzpicture}[scale = \diagramscale, vcenter=5]
\tikzgrid{5}
\filldraw[functor fill, functor draw, double] (3, 4) -- (4, 4) -- (4, 2) -- (6, 2) -- (6, 6) -- (3, 6) -- cycle;
\end{tikzpicture}
\;\; \longrightarrow \;\;
\begin{tikzpicture}[scale = \diagramscale, vcenter=5]
\tikzgrid{5}
\filldraw[functor fill, functor draw, double] (2, 4) -- (3, 4) -- (3, 2) -- (6, 2) -- (6, 6) -- (2, 6) -- cycle;
\filldraw[functor fill, functor draw, double] (3, 2) -- (4, 2) -- (4, 0) -- (6, 0) -- (6, 6) -- (3, 6) -- cycle;
\end{tikzpicture}
\;\; \longrightarrow \;\;
\begin{tikzpicture}[scale = \diagramscale, vcenter=5]
\tikzgrid{5}
\filldraw[functor fill, functor draw] (0, 2) -- (1, 2) -- (1, 0) -- (6, 0) -- (6, 6) -- (0, 6) -- cycle;
\filldraw[functor fill, functor draw] (2, 2) -- (3, 2) -- (3, 0) -- (6, 0) -- (6, 6) -- (2, 6) -- cycle;
\end{tikzpicture}
\;\; \longrightarrow \;\;
\begin{tikzpicture}[scale = \diagramscale, vcenter=5]
\tikzgrid{5}
\filldraw[functor fill, functor draw] (0, 0) -- (0, 4) -- (3, 4) -- (3, 2) -- (4, 2) -- (4, 0) -- cycle;
\end{tikzpicture}
\]
where the double lines indicate summands with multiplicity $2$.

\halfsection{Spread Modules}
\label{S:SpreadModules}

\subsection*{Collection}
We consider the following collection:
\[
\mathcal{C} := \{\coker(G \subseteq F) \mid G \subseteq F \subseteq K_I\}.
\] 
Functors in this collection have been studied in other works. They are exactly the  spread modules as defined in~\cite{BBH2022}. Recall that spread modules are defined using two antichains $S$ and $T$
in $I$ such that, for all $s$ in $S$ there exists $t$ in $T$ such that $s \leq t$, and for all $t$ in $T$ there exists $s$ in $S$ such that $s \leq t$. The \emph{spread} with \emph{sources} $S$ and \emph{sinks} $T$ is the subset of $I$
\[
[S, T] := \{v \in I \mid \exists s \in S, \exists t \in T, s \leq v \leq t\}.
\]
The \emph{spread module} with sources $S$ and sinks $T$ is then defined as the functor $K_{[S, T]}\colon I\to\vect$ where, for all $v \leq w$ in $I$,
\begin{align*}
K_{[S, T]}(v) & = \begin{cases}
K & \text{if } v \in [S, T], \\
0 & \text{otherwise,}
\end{cases}
&
K_{[S, T]}(v \leq w) & = \begin{cases}
\text{id}_K & \text{if } v, w \in [S, T], \\
0 & \text{otherwise.}
\end{cases}
\end{align*}

\begin{lemma}
\label{L:SpreadModules}
The following statements about $M\colon I\to\vect$  are equivalent:
\begin{enumerate}
\item $M$ is isomorphic to a spread module;
\item $M$ is isomorphic to $\text{\rm coker}(G\subseteq F)$ for some $ G\subseteq F\subseteq K_I$;
\item For all $v \leq w $ in $I$, $\dim(M(v))\leq 1$ and $M(v\leq w)$ is of maximal rank.
\end{enumerate}
\end{lemma}
\begin{proof}
$(1 \Rightarrow 2)$ 
If $M$ is isomorphic to a spread module $K_{[S, T]}$, then
there is a surjection $\alpha\colon K(S, -) \to M$,
and hence $M$ is isomorphic to
$\text{coker}(\text{ker}(\alpha)\subseteq K(S, -))$.
\smallskip

\noindent
$(2 \Rightarrow 3)$ Let $v \leq w$ be elements of $I$. Since $F$ surjects onto $\coker(G \subseteq F)$, we have $\dim(M(v)) \leq 1$. If $M(v) = 0$ or $M(w) = 0$, then $M(v \leq w)$ is of full rank. Otherwise, $v$ and $w$ are outside of the support of $G$, which is an upset, and so $M(x) \cong F(x)$ for all $v \leq x \leq w$, and so $M(v \leq w)$ is of full rank.
\smallskip

\noindent
$(3 \Rightarrow 1)$ Let $S$ be the minimal elements and $T$ the maximal elements of $\supp(M)$, and consider the spread module $K_{[S, T]}$. Note that, for $v \leq w$ in $I$, if $M(v \leq w)$ is nonzero, then it is an isomorphism of $1$-dimensional vector spaces. Thus, for any two elements $v, w$ in a connected component of $\supp(M)$, there exists an isomorphism $\varphi_{v, w} \colon M(v) \to M(w)$ that commutes with the transition maps of $M$. This isomorphism can be defined as a zigzag composition of transition maps.

Let $C_1, \ldots, C_k$ be the connected components of $\supp(M)$. For each $i$, fix an element $v_i \in C_i$, and define a natural transformation $f \colon M \to K_{[S, T]}$ by
\[
f(v) = \begin{cases}
\varphi_{v_i, v}(1) \id_K & \text{if } v \in C_i, i \in \{1, \ldots, k\}, \\
0 & \text{otherwise,}
\end{cases}
\]
We conclude by observing that $f$ is an isomorphism.
\end{proof}

\subsection*{Parameterization}
Consider the following poset equipped with the product order:
\[
\Omega := \{(F, G) \in \text{Sub}(I)^2 \mid F \leqJ G\}
\]
The functor  $\mathcal{Q} \colon \Omega^\op \to \Fun(I, \vect)$ which sends $(F, G)$ to  $\coker(G \subseteq F)$
is a natural  parameterization of the  collection $\mathcal{C}$.

\subsection*{Non-thinness}
The functor $\mathcal{Q}$ in general is not thin, as the following proposition illustrates, and hence we cannot apply Corollary~\ref{C:BettiFromDegeneracy}.

\begin{prop}
\label{P:SpreadModulesNotThin}
Suppose $I$ is a finite poset containing two incomparable elements $v$ and $w$.
Then the collection of nonzero elements in $\mathcal{C}$ is not independent and $\mathcal{Q} \colon \Omega^\op \to \Fun(I, \vect)$ is not thin.
\end{prop}
\begin{proof}
Let $F,G, H\colon I\to\vect$ be the unique functors in $\mathcal{C}$ whose supports are equal to the following spreads: $\text{supp}(F)=\{v\}$,  $\text{supp}(G)=\{w\}$, and
$\text{supp}(H)=\{v,w\}$. Since $H$ is isomorphic to $F\oplus G$, the collection of nonzero elements in $\mathcal{C}$ cannot be independent. Consequently $\mathcal{Q}$ cannot be thin (see Proposition~\ref{P:ThinIndependence}).
\end{proof}

Although $\mathcal{Q}$ is not thin, it can be used  to construct  thin functors.
Our strategy is to use 
Proposition~\ref{P:SingleSourceThinness} to look for thin restrictions of $\mathcal{Q}$ to subposets of $\Omega$.

\halfsection{Single-source spread modules}
\label{E:SingleSourceSpreadModules}

\subsection*{Collection}
We consider the following subcollection of $\mathcal{C}$ (see \ref{S:SpreadModules}):
\[\mathcal{C}_0:=\{F \in \mathcal{C} \mid F\text{ has exactly one generator}\}.\]
Note that all the functors  in $\mathcal{C}_0$ are of the form
$\text{coker}(G \subseteq K(v,-))$ for some $v$ in $I$. Thus functors in $\mathcal{C}_0$ are isomorphic to  spread modules $K_{[S,T]}$ with $|S|=1$. These spread modules are studied in \cite{BBH2022} and are called \textbf{single-source spread modules} . 

\subsection*{Parameterization}
The restriction of $\mathcal{Q}\colon \Omega^\op\to \Fun(I, \vect)$ 
(see~\ref{S:SpreadModules})
to the following subposet of  $\Omega$ is a natural  parameterization of the collection $\mathcal{C}_0$:
\[
\Omega_0:= \{(K(v,-), G) \in \text{Sub}(I)^2 \mid K(v,-) \leqJ G \text{ and $v$ in $I$}\}.
\]
Let $\mathcal{Q}_0 \colon \Omega_0^\op\to \Fun(I, \vect)$ denote this restriction.

\subsection*{Thinness}
The parameterization $\mathcal{Q}_0$ satisfies the conditions of Proposition~\ref{P:SingleSourceThinness}. Indeed, by definition, the functor $\mathcal{Q}_0(K(v,-), G) = \text{coker}(G \subseteq K(v,-))$
has at most one generator, and if the functor is nonzero, then its generator is at $v$. Moreover, we can check that if $(K(v,-), G) \not\leqJ (K(v', -), G')$ in $\Omega_0$, then $\Nat(\mathcal{Q}_0(K(v',-), G'), \mathcal{Q}_0(K(v, -), G)) = 0$. Consequently, $\mathcal{Q}_0$ is thin and we can utilize
Proposition~\ref{P:RelativeZerothBetti}.(1) to express
the  $0^\th$ $\mathcal{C}_0$-Betti diagrams of functors indexed by $I$ 
in terms of the standard $0^\th$ Betti diagrams of functors indexed by $\Omega_0$. 

In order to apply Corollary~\ref{C:BettiFromDegeneracy}
to use Koszul complexes to calculate the higher $\mathcal{C}_0$-Betti diagrams, 
we need to check that the poset $\Omega_0$ is an upper semilattice and that the functor $\mathcal{Q}_0$ sends $\langle \supp(\beta^0 \ker \eta_a) \rangle$ to $0$ for all $a$ in $\Omega_0$.

\subsection*{Upper semilattice}
Suppose $(I, \le)$ is an upper semilattice. 
Then $\Omega_0$ is also an upper semilattice.
To describe the join operation and the  parents of elements in $\Omega_0$,
it is convenient to identify  it with the following subposet
of  the product $I\times \Up(I)^\op$
\[\Omega_0\cong \{(v,U)\in I\times \textnormal{Up}(I)^\op\mid v\leq u \text{ for every $u$ in $U$}\} \]
by assigning to 
 an element $(K(v,-), G)$ in $\Omega_0$, the pair
$(v,\supp(G))$.
In this poset $(v_1,U_1)\leqJ (v_2,U_2)$ if and only if $v_1\leq v_2$ and $U_1\supseteq U_2$,
and hence the  join
of $(v_1,U_1)$ and $(v_2,U_2)$ in $\Omega_0$ is given by $(v_1\vee v_2,U_1\cap U_2)$.

\subsection*{Degeneracy}
By definition, for all $a = (v, U)$ in $\Omega_{0}$, the support of the kernel of the unit $\eta_{a} \colon K(a, -) \to \mathcal{RL} K(a, -)$ consists of the elements $b$ in $\Omega_{0}$ such that $b \geqJ a$ and $\mathcal{RL} K(a, b) = \Nat(\mathcal{Q}_0(b), \mathcal{Q}_0(a)) = 0$. Thus the elements in $\supp(\beta^0 \ker \eta_a)$, which are the minimal elements in $\supp(\ker \eta_a)$, are of the form $(u, (u \leq I))$ where $u$ is a minimal element of $U$.
Note further that $\mathcal{Q}_0(v', v' \leq I) = 0$ 
 for every $v'$ in $I$. Since the collection of elements in
$\Omega_{0}$ of the form $(v', (v' \leq I))$ is closed under joins, we conclude that $\mathcal{Q}_0$ sends every element of $\langle \supp(\beta^0 \ker \eta_{a}) \rangle$ to $0$.

\subsection*{Parents}
The poset $\Omega_0$ has the product order and Lemma~\ref{L:UpsetParents} gives the parents of upsets, so the parents of $(v, U)$ in $\Omega_{0}$ are of the form
\[
\mathcal{U}_{\Omega_{0}}(v, U) = (\mathcal{U}_I(v) \times \{U\}) \; \cup \; \{(v, U \cup \{w\}) \mid w \in (v \leq \Max (U^c))\},
\]
where $(v \leq \Max (U^c))$ is the set of maximal elements of $U^c$ bounded below by $v$.

\subsection*{Koszul complex}
For $M \colon I \to \vect$ the value of the functor
$\mathcal{R} M\colon \Omega_{0}\to\vect$ at 
  $(v, U)$ in $\Omega_{0}$ can be identified with
\[
\mathcal{R} M(v, U) = \Nat(\coker(K_U \subseteq K(v, -)), M) = \bigcap_{u \in \Min(U)} \ker M(v \leq u),
\]
where the second equality comes from
identifying every element $\varphi$ in $\Nat(\coker(K_U \subseteq K(v, -)), M)$ with its value $\varphi(v)(1)$ in $\bigcap_{u \in \Min(U)} \ker M(v \leq u)$.

The Koszul complex of $\mathcal{R} M$ at $(v, U)$ in degree $d \geq 0$ can now be identified with
\[
(\mathcal{K}_{(v,U)} \mathcal{R} M)_d = \bigoplus_{\substack{S \subseteq \mathcal{U}_I(v)\\ T \subseteq (v \leq \Max (U^c)) \\ |S| + |T| = d \\ S \text{ has lower bound}}} \bigcap_{u \in \Min(U \cup T)} \ker M({\textstyle \bigwedge_{(I \leq v)} S \leq u}),
\]
with differentials as defined in \ref{Std:KoszulComplexes}. Finally, we obtain the relative Betti diagrams of $M$ from the homology of this chain complex. For example, the relative projective resolution of $M_0 = \coker(K_{\{04, 32, 40\}} \subseteq K(00, -))$ is itself, since it is already a single-source spread module:
\[
\begin{tikzpicture}[scale = \diagramscale, vcenter=5]
\tikzgrid{5}
\filldraw[functor fill, functor draw] (0, 0) -- (0, 4) -- (3, 4) -- (3, 2) -- (4, 2) -- (4, 0) -- cycle;
\end{tikzpicture}
\quad \longrightarrow \quad
\begin{tikzpicture}[scale = \diagramscale, vcenter=5]
\tikzgrid{5}
\filldraw[functor fill, functor draw] (0, 0) -- (0, 4) -- (3, 4) -- (3, 2) -- (4, 2) -- (4, 0) -- cycle;
\end{tikzpicture}
\]

\halfsection{Lower hooks}
\label{E:LowerHooks}

\subsection*{Collection}
Let $(I, \leq)$ be a finite upper semilattice. We now consider the collection of lower hooks, originally defined in \cite{BOO2021}, which we define as the subcollection of $\mathcal{C}$ (see~\ref{S:SpreadModules})
\[
\mathcal{C}_L := \{\coker(K(w, -) \subseteq K(v, -)) \mid K(w, -) \subseteq K(v, -) \subseteq K_I\}.
\]

\subsection*{Parameterization}
Let $(J, \leqJ)$ be the poset $\{(v, w) \mid v, w \in I, v \leq w\}$ equipped with the product order: $(v,w)\leqJ (v',w')$ in $J$ if and only if $v\leq v'$ and $w\leq w'$ in $I$. We identify $J$ with the
subposet of the poset $\Omega_{0}$ (see~\ref{E:SingleSourceSpreadModules}) consisting of the pairs $(v, w\leq I)$ for $(v,w)$ in $J$. 
Via this identification, we parameterize lower hooks by the restriction of the functor $\mathcal{Q}_0 \colon \Omega_0^\op \to \Fun(I, \vect)$
(see~\ref{E:SingleSourceSpreadModules}) to $J$, and denote it by $\parameterization$.
Explicitly, for $a = (v,w)$ in $J$, the functor $\parameterization(a)$ is given by $\coker(K(w,-)\subseteq K(v,-))$.

\subsection*{Thinness}
The functor $\parameterization$ is thin by Proposition~\ref{P:ThinRestriction}, since it is a restriction of the thin functor $\mathcal{Q}_0$.

\subsection*{Upper semilattice}
The poset $J$ is a finite upper semilattice, inheriting the structure from $I$. Note however that $J$ is not a sublattice of $\Omega_0$.

\subsection*{Degeneracy}
For all $a = (v, w)$ in $J$, the sublattice $\langle \supp(\beta^0 \ker \eta_a) \rangle \subseteq J$ is sent to $0$ by 
$\parameterization$. Indeed, the support of the kernel of the unit $\eta_a \colon K(a, -) \to \mathcal{RL} K(a, -)$ consists of the elements $b$ of $J$ such that $b \geqJ a$ and $\mathcal{RL} K(a, b) = \Nat(\parameterization(a), \parameterization(b)) = 0$. This support has a single minimal element, $(w, w)$, which is therefore also the only element of $\supp(\beta^0 \ker \eta_a)$.
Thus, in this case, $\langle \supp(\beta^0 \ker \eta_a) \rangle = \{(w, w)\}$, which is sent to $0$ by $\parameterization$.

By Corollary~\ref{C:BettiFromDegeneracy}, we conclude that the Betti diagrams relative to lower hooks can be computed via Koszul complexes over $J$.

\subsection*{Parents}
First, for all $(v, w)$ in $J$, we have
\[
\mathcal{U}_J(v, w) = (\mathcal{U}_I(v) \times \{w\}) \cup (\{v\} \times (v \leq \mathcal{U}_I(w))),
\]
where $(v \leq \mathcal{U}_I(w))$ is the subset of $\mathcal{U}_I(w)$ of elements greater or equal to $v$.

\subsection*{Koszul complexes}
Next, we identify $\mathcal{R} M(v, w) = \Nat(\parameterization(v, w), M)$ with 
$\ker(M(v \leq w))$. Thus we can compute the Koszul complex of $\mathcal{R} M$ at $(v, w)$ in degree $d \geq 0$ as
\[
(\mathcal{K}_{(v, w)} \mathcal{R} M)_d = \bigoplus_{\substack{S \subseteq \mathcal{U}_I(v),\ T \subseteq (v \leq \mathcal{U}_I(w)) \\ |S| + |T| = d \\ S \text{ has lower bound}}} \ker M({\textstyle \bigwedge_{(I \leq v)} S \leq \bigwedge_{(I \leq w)} T}),
\]
with differentials as defined in \ref{Std:KoszulComplexes}. Finally, we obtain the relative Betti diagrams of $M$ from the homology of this chain complex.

\subsection*{Variation on lower hooks}
Now consider the subcollection of $\mathcal{C}$
\[
\mathcal{C}_\infty := \mathcal{C}_L \cup \{\coker(0 \subseteq K(v, -)) \mid v \in I\}.
\]
We parameterize this by the finite upper semilattice $J_\infty = \{(v, w) \mid v, w \in I \cup \{\infty\}, v \leq w\}$, where $\infty$ is greater than every element of $I$. 
The upper semilattice $J$ is a sublattice of $J_\infty$.
We extend the functor $\parameterization\colon J^\op\to \Fun(I, \vect)$ to $J_\infty$ by sending $(v, \infty)$ to $K(v, -)$ and denote it by the symbol $\parameterization_\infty$.
The conditions for Corollary~\ref{C:BettiFromDegeneracy} also hold for this extension. In addition, $\mathcal{C}_\infty$-exactness implies exactness and \emph{rank additivity} \cite[Prop.\ 4.3]{BOO2021}: given a short $\mathcal{C}_\infty$-exact sequence $0 \to M \to N \to L \to 0$ and $v \leq w$ in $I$, we have $\rank N(v \leq w) = \rank M(v \leq w) + \rank L(v \leq w)$.

There is a similar formula for Koszul complexes where we identify $\mathcal{R} M(v, \infty)$ with $M(v)$. For both of these examples
($\parameterization$ and $\parameterization_\infty$), the relative projective resolution of $M_0$ is
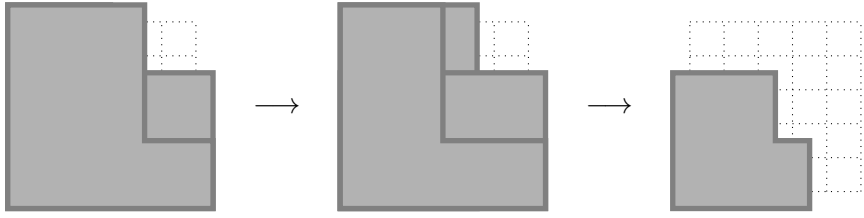
\[
\begin{tikzpicture}[scale = \diagramscale, vcenter=5]
\tikzgrid{5}
\filldraw[functor fill, functor draw] (0, 0) -- (0, 6) -- (3, 6) -- (3, 4) -- (6, 4) -- (6, 0) -- cycle;
\filldraw[functor fill, functor draw] (0, 0) -- (0, 6) -- (4, 6) -- (4, 2) -- (6, 2) -- (6, 0) -- cycle;
\end{tikzpicture}
\quad \longrightarrow \quad
\begin{tikzpicture}[scale = \diagramscale, vcenter=5]
\tikzgrid{5}
\filldraw[functor fill, functor draw] (0, 0) -- (0, 6) -- (4, 6) -- (4, 0) -- cycle;
\filldraw[functor fill, functor draw] (0, 0) -- (6, 0) -- (6, 4) -- (0, 4) -- cycle;
\filldraw[functor fill, functor draw] (0, 0) -- (0, 6) -- (3, 6) -- (3, 2) -- (6, 2) -- (6, 0) -- cycle;
\end{tikzpicture}
\quad \longrightarrow \quad
\begin{tikzpicture}[scale = \diagramscale, vcenter=5]
\tikzgrid{5}
\filldraw[functor fill, functor draw] (0, 0) -- (0, 4) -- (3, 4) -- (3, 2) -- (4, 2) -- (4, 0) -- cycle;
\end{tikzpicture}
\]
As observed above, this sequence is rank additive.

\subsection*{Relative projective dimension}
By~\cite[Corollary 10.18]{tameness} we can deduce an upper bound for the projective dimension of a functor in $\text{Fun}(J,\vect)$, i.e.\  the maximal length of its minimal projective resolution. By Corollary~\ref{C:BettiFromDegeneracy}, the maximal length of a minimal projective resolution in $\text{Fun}(I, \vect)$ relative to lower hooks will have the same upper bound. Recent work~\cite{BOOS2022} proposes 
a tighter bound in the case where $I$ is a finite lattice. Here we present a proof of this result using Koszul complexes.

\begin{prop}[cf.\ {\cite[Theorem 5.18]{BOOS2022}}]
Let $I$ be a finite lattice such that every its element has at most $n$ parents.
Then the projective dimension relative to lower hooks  of every functor 
 $M \colon I \to \vect$ is bounded above by $2n - 2$.
\end{prop}

\begin{proof}
By \cite[Corollary 10.18]{tameness}, every functor $M$ in $\Fun(I, \vect)$ has projective dimension at most $n$.
For the same reason and by observing that $J$ is also a lattice where every element has at most $2n$ parents, a functor in $\text{Fun}(J, \text{vect}_K)$ has projective dimension at most $2n$. In particular, this is the case for the functor $\mathcal{R} M \colon J \to \vect$.

Let  $\pi_1,\pi_2\colon J\to I$ be the functors defined by $\pi_1(v,w)=v$ and $\pi_2(v,w)=w$. Consider the compositions $M\pi_1$ and $M\pi_2$ in $\text{Fun}(J,\text{vect}_K)$ and the natural transformation $\alpha : M\pi_1 \to M\pi_2$ whose component at any $(v,w)$ in $J$ is $M(v\le w)$.
Observe that there is an exact sequence of functors $0 \to \mathcal{R} M = \ker(\alpha) \hookrightarrow M \pi_1 \xrightarrow{\alpha} M \pi_2$.
Now let $C\to M$ be a minimal free resolution in $\text{Fun}(I,\text{vect}_K)$.  For $i=1,2$, the sequence $C\pi_i\to M\pi_i$ in $\text{Fun}(J,\text{vect}_K)$ is also exact because exactness is defined componentwise. 
Moreover, for $u$ in $I$ and $(v, w)$ in $J$, the vector space $(K(u,-)\pi_1)(v, w)$ is nonzero if, and only if, $u \le v$, or equivalently $(u,u)\le (v,w)$. Thus, $K(u,-)\pi_1$ is isomorphic to the free functor $K((u,u), -)$. As a consequence, the sequence $C\pi_1\to M\pi_1$ is a free resolution in $\text{Fun}(J,\text{vect}_K)$, and so $M\pi_1$ has projective dimension at most $n$.
Analogously, $K(u,-)\pi_2$ coincides with the free functor $K((e,u),-)$, where $e$ is the global minimum of the lattice (it exists because $I$ is finite). We have again that the sequence $C\pi_2\to M\pi_2$ is a free resolution in $\text{Fun}(J,\text{vect}_K)$, and $M\pi_2$ has projective dimension at most $n$. 

Now consider the  exact sequence $0\to\mathcal{R}M\to M\pi_1\to \text{im}(\alpha)\to 0$ and let $(v, w)$ be an element in $J$. Since the functor $\mathcal{K}_{(v, w)}$ is exact, the sequence $0\to\mathcal{K}_{(v,w)}\mathcal{R}M\to \mathcal{K}_{(v,w)}M\pi_1 \to \mathcal{K}_{(v,w)} \text{im}(\alpha)\to 0$ is also  exact.
As we observed, the elements in the poset $J$ have at most $2n$ parents, so we obtain the following long exact sequence:
\[
\begin{tikzcd}[column sep = 3ex]
0
  \arrow[r]
&
H_{2n}(\mathcal{K}_{(v,w)}\mathcal{R}M)
  \arrow[r]
&
H_{2n}(\mathcal{K}_{(v,w)}M\pi_1)
  \arrow[r]
&
H_{2n}(\mathcal{K}_{(v,w)}\text{im}(\alpha))
  \arrow[dll, phantom, ""{coordinate, name=Z}]
  \arrow[dll, rounded corners,
    to path={ -- ([xshift=2ex]\tikztostart.east)
    |- (Z)
    -| ([xshift=-2ex]\tikztotarget.west)
    -- (\tikztotarget)}]
\\
&
H_{2n - 1}(\mathcal{K}_{(v,w)}\mathcal{R}M)
  \arrow[r]
&
H_{2n - 1}(\mathcal{K}_{(v,w)}M\pi_1)
  \arrow[r]
&
H_{2n - 1}(\mathcal{K}_{(v,w)}\text{im}(\alpha))
  \arrow[r]
&
\cdots
\end{tikzcd}
\]
However, as seen above the $2n^\th$ Betti diagram $\beta^{2n}(M\pi_1) \colon J \to \mathbb{N}$ is identically $0$, and hence $H_{2n}(\mathcal{K}_{(v,w)}M\pi_1)=0$. We deduce the equalities $\beta^{2n}\mathcal{R}M(v,w)=\dim H_{2n}(\mathcal{K}_{(v,w)}\mathcal{R}M)=0$.
Similarly, consider the short exact sequence $0\to \text{im}(\alpha)\to M\pi_2\to Q \to 0$ where $Q:=M\pi_2/\text{im}(\alpha)$. As before, it leads to a long exact sequence
\[
0
\longrightarrow H_{2n}(\mathcal{K}_{(v,w)}\text{im}(\alpha))
\longrightarrow  H_{2n}(\mathcal{K}_{(v,w)}M\pi_2)
\longrightarrow H_{2n}(\mathcal{K}_{(v,w)}Q)
\longrightarrow \cdots
\]
Consequently, $H_{2n}(\mathcal{K}_{(v,w)}M\pi_2)=0$, and so $H_{2n}(\mathcal{K}_{(v,w)}\text{im}(\alpha))=0$. Moreover, considering the first long exact sequence, we also have $H_{2n - 1}(\mathcal{K}_{(v, w)} \pi_1) = 0$, and so  $\beta^{2n-1}\mathcal{R}M(v,w)=\dim H_{2n-1}(\mathcal{K}_{(v,w)}\mathcal{R}M)=0$.

The functor $\mathcal{R} M$ therefore has projective dimension at most $2n - 2$, and so $M$ has projective dimension relative to lower hooks at most $2n - 2$.
\end{proof}

\halfsection{Rectangles}
\label{E:Rectangles}
We now illustrate an example where the assumptions of Corollary~\ref{C:BettiFromDegeneracy} fail. 

\subsection*{Collection}
Let $(I, \leq)$ be a finite upper semilattice. We define the \emph{rectangle from $v$ to $w$} as the spread $[v, w] := \{x \in I \mid v \leq x \leq w\}$ of $I$ (see~\ref{S:SpreadModules}) and we consider the subcollection of $\mathcal{C}$ of rectangle functors:
\[
\mathcal{C}_R := \{K_{[v, w]} \mid v, w \in I,\ v \leq w\}.
\]

\subsection*{Parameterization}
Let $(J, \leqJ)$ be the poset $\{(v, w) \mid v, w \in I, v \leq w\}$ with the product order, as in the example of lower hooks (see~\ref{E:LowerHooks}). This time, we identify $J$ with the subposet $\{(v, \ker(K(v, -) \to K(-, w)^*) \mid v \leq w\}$ of $\Omega_0$, where the functor $K(-, w)^*$ is equal to $K$ at $u$, when $u \leq w$, and $0$ otherwise, with transition functions of full rank. Via this identification, we parameterize rectangle functors by the restriction of the functor $\mathcal{Q}_0 \colon \Omega_0^\op \to \Fun(I, \vect)$ (see~\ref{E:SingleSourceSpreadModules}) to $J$, and denote it by $\parameterization$.

\subsection*{Thinness}
The functor $\parameterization$ is thin by Proposition~\ref{P:ThinRestriction}, since it is a restriction of the thin functor $\mathcal{Q}_0$.

\subsection*{Upper semilattice}
The poset $J$ is a finite upper semilattice just as in the previous example.

\subsection*{Problem with degeneracy locus}
Let $I$ be the finite grid $\{0 < \cdots < 5\}^2$ with the product order and $M_0\colon I\to \vect$ a functor. 
Choose $a := (04, 24)$ in $J$. It has has two parents, $(03, 24)$ and $(04, 14)$.
The Koszul complex of $\mathcal{R} M_0$ at $a$ is therefore
\[
\Nat(\parameterization(03, 14), M_0) \to
\Nat(\parameterization(03, 24), M_0) \oplus \Nat(\parameterization(04, 14), M_0) \to
\Nat(\parameterization(04, 24), M_0),
\]
which simplifies to $0 \to \Nat(\parameterization(03, 24), M_0) = K \to 0$, and so $\beta^1 \mathcal{R} M_0 (a) = 1$. 
However, in \ref{E:GridRectangles} below we compute, using a different parameterization of the collection $\mathcal{C}_R$, the relative Betti diagrams of $M_0$. In particular, we determine  $\beta^1_\mathcal{P} M_0 (\parameterization(a)) = 0$, hence $\beta^1_\mathcal{P} M_0 (\parameterization(a)) \neq \beta^1 \mathcal{R} M_0 (a)$.
The reason for this non-equality is that $\parameterization$ does not send the sublattice $\langle \supp(\beta^0 \ker \eta_a) \rangle$ to $0$ for all $a$ in $J$. In fact, $\parameterization$ does not send any element of $J$ to $0$, but sometimes $\ker \eta_a$ is nonzero, and so the support of $\beta^0 \ker \eta_a$ must be sent to nonzero functors.

\halfsection{Rectangles on a grid}
\label{E:GridRectangles}
We can amend the example of rectangles (see~\ref{E:Rectangles}) when we have more control over the poset $I$. In this example, we consider $I := \{0 < \cdots < n\}^r$ a finite grid with the product order, denoted by $\leq$. For $v, w$ in $I$, we write $v + (w - v)_i$ for the element $(v_1, \ldots, v_{i - 1}, w_i, v_{i + 1}, \ldots, v_r)$.

\subsection*{Collection}
We consider the same collection $\mathcal{C}_R$ as in~\ref{E:Rectangles}.

\subsection*{Parameterization}
Let $(J, \leqJ)$ be the poset $\{(v, w) \mid v, w \in I, v \leq w\}$ with the product order, as in the previous example. This time, we identify $J$ with the subposet 
\[
\Big\{\Big(v, \bigcup_{i = 1}^r(v + (w - v)_i \leq I)\Big) \;\Big|\; v, w \in I, v \leq w\Big\}
\]
of $\Omega_0$ (see~\ref{E:SingleSourceSpreadModules}). 
Via this identification, we parameterize rectangle functors by the restriction of the functor $\parameterization \colon \Omega_0^\op \to \Fun(I, \vect)$ (see~\ref{E:SingleSourceSpreadModules}) to $J$, and denote it by $\parameterization$. Explicitly, $\parameterization$ maps $(v, w)$ to the functor $\coker(\bigoplus_{i = 1}^r K(v + (w - v)_i, -) \rightarrow K(v, -))$. 

\subsection*{Thinness}
The functor $\parameterization$ is thin by Proposition~\ref{P:ThinRestriction}, since it is a restriction of the thin functor $\mathcal{Q}_0$.

\subsection*{Upper semilattice}
The poset $J$ is a finite upper semilattice just as in the previous example.

\subsection*{Degeneracy}
For all $a = (v, w)$ in $J$, $\parameterization$ sends the sublattice $\langle \supp(\beta^0 \ker \eta_a) \rangle$ to $0$. Indeed, the support of the kernel of the unit $\eta_a \colon K(a, -) \to \mathcal{RL} K(a, -)$ consists of the elements $b$ in $J$ such that $b \geqJ a$ and $\mathcal{RL} K(a, b) = \Nat(\parameterization(b), \parameterization(a)) = 0$. Thus the elements of $\supp(\beta^0 \ker \eta_a)$, which are the minimal elements of $\supp(\ker \eta_a)$, are of the form $(v + (w - v)_i, w)$ where $i \in \{1, \ldots, r\}$.
The sublattice $\langle \supp(\beta^0 \ker \eta_a) \rangle$ then consists of the elements $(v_X , w)$, for every nonempty subset $X$ of $\{1, \ldots, r\}$, where, for all $i$ in $\{1, \ldots, r\}$, the element $v_X$ in $I$ is defined by
\[
(v_X)_i := \begin{cases}
w_i & \text{if } i \in X, \\
v_i & \text{otherwise.}
\end{cases}
\]
Note that for $(v', w')$ in $J$, $\parameterization(v', w') = 0$ if, and only if, there exists $i \in \{1, \ldots, r\}$ such that $v'_i = w'_i$. As a consequence, $\parameterization$ sends the elements of the sublattice $\langle \supp(\beta^0 \ker \eta_a) \rangle$ to $0$.

By Corollary~\ref{C:BettiFromDegeneracy}, we can compute Betti diagrams relative to rectangles via Koszul complexes.

\subsection*{Parents}
Parents are exactly the same as for lower hooks (see~\ref{E:LowerHooks}).

\subsection*{Koszul complexes}
For $(v, w)$ in $J$, we can identify $\mathcal{R} M(v, w) = \Nat(\parameterization(v, w), M)$ with $\bigcap_{i = 1}^r \ker M(v \leq v + (w - v)_i)$. Thus we can compute the Koszul complex of $\mathcal{R} M$ at $(v, w)$ in degree $d \geq 0$ as
\[
(\mathcal{K}_{(v, w)} \mathcal{R} M)_d = \bigoplus_{\substack{S \subseteq \mathcal{U}_I(v) \\ T \subseteq (v \leq \mathcal{U}_I(w)) \\ |S| + |T| = d }} \bigcap_{i = 1}^r \ker M\left(v_S \leq v_S + (w_T - v_S)_i\right),
\]
where $v_S := \bigwedge_{(I \leq v)} S$ and $w_T := \bigwedge_{(I \leq w)} T$, and
with differentials as defined in \ref{Std:KoszulComplexes}. We obtain the relative Betti diagrams of $M$ from the homology of this chain complex. For example, the relative projective resolution of $M_0$ is
\[
\begin{tikzpicture}[scale = \diagramscale, vcenter=5]
\tikzgrid{5}
\filldraw[functor fill, functor draw, line cap = round] (3, 2) rectangle (4, 4);
\end{tikzpicture}
\quad \longrightarrow \quad
\begin{tikzpicture}[scale = \diagramscale, vcenter=5]
\tikzgrid{5}
\filldraw[functor fill, functor draw] (0, 0) rectangle (4, 4);
\end{tikzpicture}
\quad \longrightarrow \quad
\begin{tikzpicture}[scale = \diagramscale, vcenter=5]
\tikzgrid{5}
\filldraw[functor fill, functor draw] (0, 0) -- (0, 4) -- (3, 4) -- (3, 2) -- (4, 2) -- (4, 0) -- cycle;
\end{tikzpicture}
\]

\bibliographystyle{plain}
\bibliography{bibliography.bib}
\end{document}